\RequirePackage[english]{babel}
\documentclass[a4paper,twoside,11pt]{amsart}
\usepackage[latin1]{inputenc}
\usepackage{amsmath}
\usepackage{amssymb}
\usepackage{amsfonts}
\usepackage{mathrsfs}
\usepackage{fixltx2e}[2005/12/01]
\usepackage{hyperref}

\usepackage{enumerate}
\renewcommand{\theenumi}{{\upshape{(\roman{enumi})}}}

\newcommand{\myitem}[1]{\renewcommand{\theenumi}{{\upshape{(#1)}}}\item}

\usepackage[matrix,arrow,cmtip,frame,graph,2cell]{xy} 
\SelectTips{cm}{}
\newdir{(}{{}*!/-5pt/@^{(}}
\newdir{(x}{{}*!/-5pt/@_{(}}
\newdir{+}{{}*!/-9pt/{}}
\newdir{>+}{@{>}*!/-9pt/{}}
\entrymodifiers={+!!<0pt,\fontdimen22\textfont2>}
\UseTwocells

\usepackage[nosubsections]{mytheorems2}

\newcommand\inj{\hookrightarrow}
\newcommand\surj{\twoheadrightarrow}
\newcommand\iso{\cong} 
\newcommand\pref[1]{\eqref{#1}}
\newcommand\etale{\'{e}tale}
\DeclareMathOperator{\Hom}{Hom}
\DeclareMathOperator{\Sym}{Sym}
\newcommand\map[3]{#1\colon #2\rightarrow #3}
\newcommand\injmap[3]{#1\colon #2\hookrightarrow #3}
\newcommand\id[1]{\mathrm{id}_{#1}} 
\newcommand\cons{\mathrm{cons}} 
\newcommand\red{\mathrm{red}}
\newcommand\sA{\mathcal{A}}
\newcommand\sO{\mathcal{O}}
\newcommand\sL{\mathcal{L}}
\newcommand\Spec{\mathrm{Spec}}
\newcommand\Proj{\mathrm{Proj}}
\newcommand\Hilb{\mathrm{Hilb}}
\newcommand\A[1]{\mathbb{A}^{#1}}    
\newcommand\SG[1]{{\mathfrak{S}_{#1}}}   
\newcommand\Z{\mathbb{Z}}
\newcommand\Q{\mathbb{Q}}

\newcommand{\fpr}{\mathrm{fpr}} 
\newcommand{\stab}{\mathrm{stab}} 
\newcommand{\metale}{\text{\'et}} 
\newcommand{\stX}{\mathscr{X}}
\newcommand{\stY}{\mathscr{Y}}
\newcommand{\stZ}{\mathscr{Z}}
\newcommand{\stW}{\mathscr{W}}

\newcommand{\stH}{\mathscr{H}}
\newcommand{\stQ}{\mathscr{Q}}

\newcommand{\cR}{\mathcal{R}} 
\newcommand{\Isom}{\mathcal{I}som}
\newcommand{\GC}{{\upshape GC}}
\newcommand{\AF}{{\upshape AF}}
\newcommand{\norm}{\mathrm{N}}
\newcommand{\trace}{\mathrm{Tr}} 
\newcommand\Sect{\mathrm{\Pi}}
\newcommand\weilr{\mathbf{R}} 

\newcommand\loccit{\textit{loc.\ cit.}}
\newcommand\catC{\mathbf{C}}
\newcommand\Et{\mathbf{\acute{E}t}}
\newcommand\Cat{\mathbf{Cat}}
\newcommand\AlgSp{\mathbf{AlgSp}}
\newcommand{\all}{\mathrm{all}}
\newcommand{\qc}{\mathrm{qc}}
\newcommand{\qsep}{\mathrm{qsep}}
\newcommand{\sep}{\mathrm{sep}}
\newcommand{\finite}{\mathrm{fin}}

\newcommand{\op}{\mathrm{op}}

\newcommand{\groupoidMAPS}[4]{\xymatrix@1@M=0mm@C=10mm{#1%
\ar@<.5ex>@{+->+}[r]^-{#3} \ar@<-.5ex>@{+->+}[r]_-{#4} & #2}}
\newcommand{\groupoid}[2]{\xymatrix@1@M=0mm@C=7mm{#1
\ar@<.5ex>@{+->+}[r] \ar@<-.5ex>@{+->+}[r] & #2}}

\begin{document}

\title[Geometric quotients]{Existence and properties\\of geometric quotients}
\author{David Rydh}
\address{Department of Mathematics, KTH, 100 44 Stockholm, Sweden}
\email{dary@math.kth.se}
\date{2012-05-04}


\begin{abstract}
In this paper, we study quotients of groupoids and coarse moduli spaces of
stacks in a general setting. Geometric quotients are not always categorical,
but we present a natural topological condition under which a geometric quotient
is categorical.
%
We also show the existence of geometric quotients of finite flat groupoids
and give explicit local descriptions.
Exploiting similar methods, we give an easy proof of the existence of quotients
of flat groupoids with finite stabilizers.
As the proofs do not use noetherian methods and are valid for general algebraic
spaces and algebraic stacks, we obtain a slightly improved version of Keel
and Mori's theorem.
\end{abstract}

\maketitle


\setcounter{secnumdepth}{0}
\begin{section}{Introduction}
Quotients appear frequently in almost every branch of algebraic geometry, most
notably in moduli problems. It has been known for a long time that the category
of algebraic spaces is better suited for quotient problems than the category of
schemes. In fact, even the easiest instances of quotients, such as quotients of
schemes by free actions of finite groups, need not exist as schemes. In the
category of algebraic spaces, on the other hand, quotients of free actions exist
almost by definition.

In the category of schemes, it is easy to prove that geometric quotients are
categorical and hence unique~\cite{GIT}. Surprisingly, this is not the case in
the category of algebraic spaces as pointed out by
Koll\'ar~\cite{kollar_quotients}. In particular, geometric quotients need not
be unique.

The purpose of this paper is threefold. Firstly, we introduce a natural
subclass of geometric quotients which we call \emph{strongly geometric} and
show that these quotients are categorical among algebraic spaces. Secondly, we
show the existence of quotients by finite groups and give explicit
\etale{}-local descriptions. Thirdly, we show the existence of quotients of
arbitrary flat groupoids with finite stabilizers, generalizing Keel and Mori's
theorem~\cite{keel_mori_quotients}.

A geometric quotient $X\to X/G$ is required to be \emph{topological}, that is,
the fibers should be the orbits and the quotient should have the quotient
topology. A strongly geometric quotient is a geometric quotient that is
\emph{strongly topological}. This means that in addition to requiring that the
quotient has the correct topology we also require that the equivalence
relation determined by the quotient has the correct topology.

Deligne has proved the existence of geometric quotients of separated algebraic
spaces by arbitrary actions of finite groups, but without any
published proof, cf.~\cite[p.~183]{knutson_alg_spaces}. Deligne's idea was to
use \emph{fixed-point reflecting} \etale{} covers to deduce the existence from
the affine case. Matsuura independently proved the aforementioned result
using the same method~\cite{matsuura_quotients}.
Koll\'ar developed Deligne's ideas in~\cite{kollar_quotients} and showed that a
geometric quotient of a proper group action is categorical in the category of
algebraic spaces.

We extend Deligne and Koll\'ar's results on fixed-point reflecting coverings.
Using effective descent results for
submersions~\cite{rydh_submersion_and_descent} we then show that strongly
geometric quotients are categorical in almost every situation\
(Theorem~\ref{T:(strong)-geom+univ-open-is-(strong)-cat}). We also settle
Koll\'ar's conjecture~\cite[Rmk.~2.20]{kollar_quotients} that
geometric quotients are categorical among \emph{locally separated} algebraic
spaces.

We then proceed to show the existence of quotients of separated algebraic
spaces by finite groups. Important to us is that we obtain an explicit
\etale{}-local description of the quotient space which does not follow from
Keel and Mori's result. This local description is particularly nice for the
symmetric product\ (Corollary~\ref{C:existence-of-Sym:sep-alg-space}) and
such an explicit description is needed to deduce properties of $\Sym^n(X/S)$
from the affine case, as is done in~\cite{skjelnes_ekedahl_good_component,
rydh_famzerocycles-I,rydh-skjelnes_gen_etale_fam}. 

The techniques developed in this paper lead to a simple proof of Keel and
Mori's theorem (Theorem~\ref{T:KM-general}).

\begin{theorem*}
Let $S$ be an algebraic space and let $\groupoidMAPS{R}{X}{s,t}{}$ be a groupoid
of algebraic spaces over $S$ such that $s$ and $t$ are flat and locally of
finite presentation and let $\map{j=(s,t)}{R}{X\times_S X}$.
If the stabilizer
$j^{-1}\bigl(\Delta(X)\bigr)\to X$ is \emph{finite}, then there is a uniform
strongly geometric and categorical quotient $X\to X/R$ such that:
%
\begin{enumerate}
\item $X/R\to S$ is separated (resp.\ quasi-separated) if and only if $j$ is
finite (resp.\ quasi-compact);
\item $X/R\to S$ is locally of finite type if $S$ is locally noetherian
and $X\to S$ is locally of finite type; and
\item $R\to X\times_{X/R} X$ is proper.
\end{enumerate}
\end{theorem*}

Note that we are able to remove the finiteness hypotheses present in earlier
proofs.
In the original theorem~\cite[Thm.~1]{keel_mori_quotients}, it is assumed
that the base scheme $S$ is locally noetherian, that $X\to S$ is locally of
finite type and that $j$ is quasi-compact. In Conrad's
treatment~\cite{conrad_coarsespace}, the base scheme $S$ is arbitrary, but
$X\to S$ is locally of finite presentation and $j$ is quasi-compact.

The hypothesis that the stabilizer is finite implies that the diagonal $j$ is
separated and locally quasi-finite~\cite[Lem.~2.7]{keel_mori_quotients}. The
stack
$\stX=[\groupoid{R}{X}]$ is thus an algebraic stack with separated and
locally quasi-finite
diagonal. The quotient $X/R$ is the \emph{coarse moduli space} of $\stX$. The
stabilizer is a pull-back of the \emph{inertia stack} $I_{\stX}\to \stX$.
Rephrased in the language of stacks our generalization of the Keel--Mori theorem
takes the following form.

\begin{theorem*}
Let $\stX$ be an algebraic stack. A coarse moduli space $\map{\pi}{\stX}{X}$
such that $\pi$ is separated exists if and only if $\stX$ has finite inertia.
In particular, any separated Deligne--Mumford stack has a coarse moduli space.
\end{theorem*}


This paper resulted from an attempt to understand fundamental questions about
group actions and quotients. As a consequence,
Sections~\ref{S:groupoids}--\ref{S:fpr-and-descent} are written in a more
general setting than needed for our generalization of the Keel--Mori theorem.
For example, the results apply to non-flat equivalence
relations~\cite{kollar_finite-equiv-rels}.
The impatient reader mainly interested in the Keel--Mori theorem is encouraged
to go directly to~\S\ref{S:cms-stacks}.

The first step in the proof of the Keel--Mori theorem is to find a flat and
locally quasi-finite presentation of a stack $\stX$ with quasi-finite
diagonal. This result is classical for noetherian stacks,
cf.\ \cite[Lem.~7.2]{gabriel_quotients}
and~\cite[Lem.~3.3]{keel_mori_quotients}. The general case can be found
in~\cite[Thm.~7.1]{rydh_etale-devissage} and~\cite[Lem.~04N0]{stacks-project}.

\vspace{-1 mm}
\begin{subsection}{Assumptions and terminology}
We work with general algebraic spaces and algebraic stacks as
in~\cite{stacks-project,rydh_etale-devissage,rydh_hilbert}. In particular, we
do not assume, as in~\cite{knutson_alg_spaces,laumon}, that algebraic spaces
are quasi-separated and that algebraic stacks have quasi-compact and separated
diagonals.
If the inertia stack is finite, then the diagonal is locally quasi-finite and
separated but not necessarily quasi-compact. Nevertheless, the reader may, if
so inclined, assume that all algebraic spaces and stacks have quasi-compact and
separated diagonals. This would, however, necessitate tedious verifications
that certain equivalence relations and groupoids have quasi-compact
diagonals, e.g., in Theorem~\pref{T:existence-of-GC-by-fpr}. The
quasi-separatedness is typically better treated at subsequent stages using
global methods.

We follow the terminology of EGA with one exception. As
in~\cite{raynaud_hensel_rings,laumon} we mean by unramified a
morphism locally of finite type and formally unramified but not necessarily
locally of finite presentation.
By a \emph{presentation} of an algebraic stack $\stX$, we always mean an
algebraic space $U$ together with a faithfully flat
morphism $U\to \stX$ locally of finite presentation.
\end{subsection}

\begin{subsection}{Structure of the article}
We begin with some general definitions and properties of quotients in
\S\S\ref{S:groupoids}\nobreakdash--\ref{S:quotients}. Quotients are treated in full
generality and we do not assume that the groupoids are flat.
From Section~\ref{S:affine-case} and on, all groupoids are flat and locally of
finite presentation.

In \S\ref{S:fpr-and-descent} we generalize the results of
Koll\'ar~\cite{kollar_quotients} on topological quotients, fixed-point
reflecting
morphisms and the \emph{descent condition}. We show that universally open
strongly geometric quotients satisfy the descent condition and hence are
categorical.
An important observation is that the fixed-point reflecting locus of a
square \etale{}
morphism between groupoids is open if the stabilizer of the target groupoid is
proper. This is the only place where the properness of the stabilizer is
utilized.

In \S\ref{S:affine-case} we give an overview of well-known results on the
existence and properties of quotients of affine schemes by finite flat
groupoids. In this generality, the results are due to Grothendieck.

In \S\ref{S:finite-quotients} we use the results of \S\ref{S:fpr-and-descent}
to deduce the existence of finite quotients for arbitrary algebraic spaces from
the affine case treated in \S\ref{S:affine-case}. The fixed-point reflecting
\etale{} cover with an essentially affine scheme that we need is constructed
as a Weil restriction.

In \S\ref{S:cms-stacks} we restate the results of \S\ref{S:fpr-and-descent} in
terms of stacks. We then deduce the existence of a coarse moduli space to any
stack $\stX$ with finite inertia stack, from the case where $\stX$ has a finite
flat presentation. Here, what we need is a fixed-point reflecting \etale{}
cover
$\stW$ of the stack $\stX$ such that the cover $\stW$ admits a finite flat
presentation. After finding a quasi-finite flat presentation, this is
accomplished using Hilbert schemes.
We also give an example of a stack with quasi-finite, but not finite, inertia
that does not have a coarse moduli space.

The existence results of \S\ref{S:finite-quotients} follow from the
independent and more general results of \S\ref{S:cms-stacks} but the
presentations of the quotients are quite different.
In~\S\ref{S:finite-quotients} we begin with an algebraic space $X$ with an
action of a finite groupoid and construct an essentially affine cover $U$ with
an action of the same groupoid. In~\S\ref{S:cms-stacks} we begin by modifying
the groupoid obtaining a quasi-finite groupoid action on an essentially affine
scheme $X$. Then we take a covering which has an action of a finite groupoid.

In the appendices we have collected some results on effective descent of
\etale{} morphisms, AF-schemes, strong
homeomorphisms and Weil restrictions.
\end{subsection}

\begin{subsection}{Acknowledgment}
I would like to thank J.\ Alper, O.\ Gabber, J.\ Hall, Y.\ Huang, D.\ Laksov,
Y.\ Matsuura, R.\ Skjelnes and two referees for many useful suggestions
and comments.
\end{subsection}

\end{section}
\setcounter{secnumdepth}{3}


\begin{section}{Groupoids and stacks}
\label{S:groupoids}

Let $G$ be a group scheme over $S$, or more generally an algebraic group space
(a group object in the category of algebraic spaces), and let $X$ be an
algebraic space over $S$. An action of $G$ on $X$ is a morphism
$\map{\sigma}{G\times_S X}{X}$
compatible with the group structure on $G$.
The group action $\sigma$ gives rise to a pre-equivalence relation
$\groupoidMAPS{G\times_S X}{X}{\sigma}{\pi_2}$, where $\pi_2$ is the second
projection, i.e., a \emph{groupoid in algebraic spaces}.

\begin{definition}
Let $S$ be an algebraic space. An \emph{$S$-groupoid in algebraic spaces}
consists of two algebraic $S$-spaces, $R$ and $U$, together with
morphisms
\begin{enumerate}
\item source and target $\map{s,t}{R}{U}$,
\item composition $\map{c}{R\times_{t,U,s} R}{R}$,
\item identity $\map{e}{U}{R}$ and
\item inverse $\map{i}{R}{R}$,
\end{enumerate}
such that $\bigl(R(T),U(T),s,t,c,e,i\bigr)$ is a groupoid in sets for every
affine $S$-scheme $T$ in a functorial way. We will denote the
groupoid by
$\groupoidMAPS{R}{U}{s}{t}$ or $(R,U)$. A morphism of groupoids
$\map{f}{(R,U)}{(R',U')}$ consists of two morphisms $R\to R'$ and $U\to U'$,
also denoted $f$, such that $\map{f}{(R(T),U(T))}{(R'(T),U'(T))}$ is a
morphism of groupoids in sets for every affine $S$-scheme~$T$.
\end{definition}

\begin{remark}
The inverse $\map{i}{R}{R}$ is an involution such that $s=t\circ i$.
Thus $s$ has a property if and only if $t$ has the same property.
Let $G\to S$ be a group scheme acting on an algebraic space $X\to S$ and let
$\groupoidMAPS{G\times_S X}{X}{s}{t}$ be the associated groupoid. If $G\to S$
has a
property stable under base change, then $s$ and $t$ have the same property.
\end{remark}

\begin{notation}\label{N:stabilizer}
By a \emph{groupoid} we will always mean a groupoid in algebraic spaces. If
$\groupoid{R}{U}$ is a groupoid, then we let $j$ be the diagonal
morphism
$$\map{j=(s,t)}{R}{U\times_S U}.$$
Let $T$ be an $S$-space and let $\map{g}{T}{U}$ be a morphism.  The
\emph{stabilizer} of $g$, denoted $\stab(g)$, is the pull-back of $j$ along
$\map{(g,g)}{T}{U\times_S U}$ and is a group space over $T$. If
$\map{h}{T'}{T}$,
then $\stab(g\circ h)=\stab(g)\times_T T'$. The \emph{stabilizer of the
groupoid} is the universal stabilizer
$${\map{\stab(U):=\stab(\id{U})}{j^{-1}\bigl(\Delta(U)\bigr)}{U}}.$$
If
$\map{r}{T}{R}$
is a morphism, then conjugation by $r$ induces an isomorphism
$\stab(s\circ r)\to \stab(t\circ r)$. An inverse is given by conjugation by
$i\circ r$.

We say that
$\groupoid{R}{U}$ is flat, locally of finite presentation,
quasi-finite, etc.\ if $s$, or equivalently $t$, is flat, locally of finite
presentation, quasi-finite, etc.
\end{notation}

\begin{xpar}[Topological equivalence relation induced by a groupoid]
Let $\groupoid{R}{U}$ be a groupoid and let $u_1,u_2\in |U|$ be points.
We say that $u_1$ and $u_2$ are \emph{$R$-equivalent}, written $u_1\sim_R u_2$,
if there exists a point $r\in
|R|$ such that $s(r)=u_1$ and $t(r)=u_2$. It follows immediately from the
groupoid axioms that $\sim_R$ is an equivalence relation on $|U|$ since
$|R\times_U
R|\to |R|\times_{|U|} |R|$ is surjective. The orbit $R(u)$ of $u\in |U|$
is the $R$-equivalence class of $u$. Explicitly we have that
$R(u)=t(s^{-1}(u))$. The map of sets $|U|\to |U|/\sim_R$ is the coequalizer of
$\groupoid{|R|}{|U|}$.
We say that a subset $W\subseteq |U|$ is \emph{$R$-stable} if $W$ is closed
under $\sim_R$, or equivalently, if $s^{-1}(W)=t^{-1}(W)$ as subsets of $|R|$.
\end{xpar}

\begin{remark}
If $\groupoid{R}{U}$ is a groupoid that is flat and locally of
finite presentation, then the associated fppf stack
$[\groupoid{R}{U}]$ is algebraic~\cite[Cor.~10.6]{laumon}
or~\cite[Thm.~06DC]{stacks-project}.
In particular, when $G\to S$ is a flat group scheme of finite presentation
acting on an algebraic space
$X\to S$, the stack
$[X/G]=[\groupoid{G\times_S X}{X}]$ is algebraic.
\end{remark}

\begin{definition}
Let $\stX$ be an algebraic stack. The \emph{inertia stack} $I_{\stX}\to \stX$ is
the pull-back of the diagonal $\Delta_{\stX/S}$ along the same morphism
$\Delta_{\stX/S}$. The inertia stack is independent of the base $S$.
\end{definition}

\begin{remark}\label{R:stacks-vs-groupoids:diagonal/inertia}
Let $\stX$ be an algebraic stack with a presentation
$\map{p}{U}{\stX}$ and let $R=U\times_\stX U$. Then
$\stX\iso[\groupoid{R}{U}]$ and we have $2$-cartesian diagrams
\begin{equation*}
\vcenter{\xymatrix@C+5mm{
R\ar[r]^-j\ar[d] & U\times_S U\ar[d]^{p\times p}\\
\stX\ar[r]^-{\Delta_{\stX/S}} & \stX\times_S \stX\ar@{}[ul]|\square
}} \quad\text{and}\quad
\vcenter{\xymatrix{
\stab(U)\ar[r]\ar[d] & U\ar[d]^p\\
I_\stX\ar[r] & \stX.\ar@{}[ul]|\square
}}
\end{equation*}
\end{remark}

\end{section}


\begin{section}{General remarks on quotients}
\label{S:quotients}

\begin{subsection}{Topological, geometric and categorical quotients}
Recall that a morphism of topological spaces $\map{f}{X}{Y}$ is
\emph{submersive} if it is surjective and $Y$ has the quotient topology, i.e., a
subset $Z\subseteq Y$ is open if and only if its inverse image $f^{-1}(Z)$ is
open. Equivalently, $Z\subseteq Y$ is closed if and only if $f^{-1}(Z)$ is
closed.

\begin{xpar}[Constructible topology]\label{X:cons-top}
Let $\map{f}{X}{Y}$ be a morphism of algebraic spaces. We say that $f$ is
submersive if the associated map of topological spaces $\map{|f|}{|X|}{|Y|}$ is
submersive~\cite[15.7.8]{egaIV}. By slight abuse of notation we say that
$f^{\cons}$ is submersive if the induced map of topological spaces
$\map{|f|^{\cons}}{|X|^{\cons}}{|Y|^{\cons}}$ is submersive. Here ``$\cons$''
denotes the \emph{constructible topology}, cf.~\cite[7.2.11]{egaI_NE}
or~\cite[1.9.13]{egaIV} and~\cite[\S1]{rydh_submersion_and_descent}. If $\stX$
is an algebraic stack (or an algebraic space) and $\map{p}{U}{\stX}$ is a
presentation, then $W\subseteq |\stX|$ is open in the constructible topology if
and only if $p^{-1}(W)$ is open in the constructible topology. A subset
$W\subseteq |\stX|$ is
constructible if and only if it is open and closed in the constructible
topology.
\end{xpar}

The groupoid-averse reader may benefit from the comparison of the following
definition with the corresponding Definition~\pref{D:quotient-of-stacks} for
stacks.

\begin{definition}
Let $\groupoid{R}{X}$ be an $S$-groupoid. A morphism
$\map{q}{X}{Y}$ is \emph{equivariant} if $q\circ s=q\circ t$. Then
$\groupoid{R}{X}$ is also a $Y$-groupoid and $\map{j}{R}{X\times_S X}$
factors through $X\times_Y X\inj X\times_S X$. We denote the morphism $R \to
X\times_Y X$ by $j_{/Y}$.

If $Y'\to Y$ is a morphism, we let $R'=R\times_Y Y'$ and $X'=X\times_Y Y'$. Then
$\map{q'}{X'}{Y'}$ is equivariant with respect to the groupoid
$\groupoid{R'}{X'}$. If
a property of $q$ is stable under flat base change $Y'\to Y$, then we say that
the
property is \emph{uniform}. If it is stable under arbitrary base change, we
say that it is \emph{universal}.
If $q$ is an equivariant morphism, then we say that
\begin{enumerate}
\item $q$ is a \emph{categorical quotient} (with respect to a full subcategory
$\catC$ of the category of algebraic spaces) if $q$ is an initial object among
equivariant morphisms (with target in $\catC$). Concretely, this means that for
any equivariant morphism $\map{r}{X}{Z}$ (with $Z\in\catC$) there is a unique
morphism $Y\to Z$ such
that the diagram
$$\xymatrix{
{X}\ar[r]^{q}\ar[rd]_{r} & {Y}\ar@{-->}[d]\\
& {Z} 
}$$
commutes.
\item $q$ is a \emph{Zariski quotient} if $\map{|q|}{|X|}{|Y|}$ is the
coequalizer of $\groupoid{|R|}{|X|}$ in the category of topological
spaces. Equivalently, the fibers of $q$ are the $R$-orbits of $|X|$ and $q$
is submersive.
\item $q$ is a \emph{constructible quotient} if
$\map{|q|^\cons}{|X|^\cons}{|Y|^\cons}$ is the coequalizer of
$\groupoid{|R|^{\cons}}{|X|^{\cons}}$ in the category of topological
spaces. Equivalently the fibers of $q$
are the $R$-orbits of $|X|$ and $q^\cons$ is submersive.
\item $q$ is a \emph{topological quotient} if it is both a universal Zariski
quotient and a universal constructible quotient.
\item $q$ is a \emph{strongly topological quotient} if it is a topological
quotient and $\map{j_{/Y}}{R}{X\times_Y X}$ is universally submersive.
\item $q$ is a \emph{geometric quotient} if it is a topological quotient
and $\sO_Y=(q_*\sO_X)^R$, i.e., if the sequence of sheaves in the \etale{}
topology
\begin{equation}\label{E:geometric-quot-seq}
\vcenter{%
\xymatrix{\sO_Y\ar[r] & q_*\sO_X\ar@<.5ex>[r]^-{s^*}\ar@<-.5ex>[r]_-{t^*} &
  (q\circ s)_*\sO_R}}
\end{equation}
is exact.
%
\item $q$ is a \emph{strongly geometric quotient} if it is geometric and
strongly topological.
\end{enumerate}
\end{definition}

\begin{lemma}
Let $\map{q}{X}{Y}$ be an equivariant morphism and let as before
$\map{j_{/Y}}{R}{X\times_Y X}$ denote the diagonal. The following are
equivalent.
\begin{enumerate}
\item\label{LI:set:orbits} For any field $k$ and point $\map{y}{\Spec(k)}{Y}$,
  the set $|X\times_Y \Spec(k)|$ has at least (resp.\ at most, resp.\ exactly)
  one $R\times_Y \Spec(k)$-orbit.
\item\label{LI:set:surjectivity} The morphism $q$ is surjective (resp.\ the
  morphism $j_{/Y}$ is surjective, resp.\ the morphisms $q$ and $j_{/Y}$ are
  surjective).
\end{enumerate}
If, in addition, $q$ and $j_{/Y}$ are locally of finite type or integral, then
these
two conditions are equivalent to:
\begin{enumerate}\setcounter{enumi}{2}
\item\label{LI:set:eqclasses} For any algebraically closed field $K$, the map
  $\map{q(K)}{X(K)/R(K)}{Y(K)}$ is surjective (resp.\ injective,
  resp.\ bijective).
\end{enumerate}
\end{lemma}
\begin{proof}
That \ref{LI:set:surjectivity}$\implies$\ref{LI:set:orbits} follows from the
definitions. We will show that
\ref{LI:set:orbits}$\implies$\ref{LI:set:surjectivity} and
\ref{LI:set:eqclasses}$\implies$\ref{LI:set:surjectivity} without any
assumptions on $q$ and $j_{/Y}$. If $q(K)$ is surjective or if every fiber has
at least one orbit, then clearly $q$ is surjective. Assume instead that $q(K)$
is injective or that every fiber has at most one orbit. Choose a point
$\map{z}{\Spec(K)}{X\times_Y X}$ and let
$\map{y}{\Spec(K)}{Y}$ be the induced point of $Y$. We also let $z$ denote the
induced point $\Spec(K)\to X_y\times_{\Spec(K)} X_y$ and let
$\map{z_1,z_2}{\Spec(K)}{X_y}$ be the induced points of $X_y$. By assumption,
there is a point $r\in |R_y|$, or even a point $r\in R_y(K)$, such that
$s(r)=z_1$ and $t(r)=z_2$ in $|X_y|$.
As the residue fields of $z_1$ and $z_2$ are $K$, we have that $z$ is the
unique point of $|X_y\times_{\Spec(K)} X_y|$ above $(z_1,z_2)\in |X_y|\times
|X_y|$. It follows that the image of $r$ along
$\map{j_y}{R_y}{X_y\times_{\Spec(K)} X_y}$ is $z$ and we conclude that $j_{/Y}$
is surjective.

Finally, assume that $q$ and $j_{/Y}$ are locally of finite type or integral
and that $K$ is an algebraically closed field. Then
\ref{LI:set:surjectivity}$\implies$\ref{LI:set:eqclasses} since the
surjectivity of $q$ and $j_{/Y}$ implies the surjectivity of $q(K)$ and
$j_{/Y}(K)$ respectively.
\end{proof}

The lemma gives the following alternative description of topological
quotients:
\begin{align*}
\text{$q$ univ.\ Zariski} &\iff  \text{$q$ univ.\ submersive and $j_{/Y}$
surjective}\\
\text{$q$ topological} &\iff \text{$q$, $q^\cons$ univ.\
submersive and $j_{/Y}$ surjective}\\
\text{$q$ strongly topological} &\iff \text{$q$, $q^\cons$, $j_{/Y}$ univ.\
submersive.}
\end{align*}

\begin{proposition}\label{P:quotients-misc-props}
Let $\map{q}{X}{Y}$ be a topological quotient of the groupoid
$\groupoid{R}{X}$.
If $s$ has one of the properties: universally
open, universally closed, quasi-compact; then so does $q$. For the first two
properties, it is sufficient that $q$ is a universal Zariski quotient.
%
\end{proposition}
\begin{proof}
If $s$ is universally open
(resp.\ universally closed, resp.\ quasi-compact), then so is the projection
$\map{\pi_1}{X\times_Y X}{X}$ as $j_{/Y}$ is surjective.
Since $q$ and $q^\cons$ are universally submersive,
we have that $q$ is universally open (resp.\ universally closed,
resp.\ quasi-compact) by~\cite[Prop.~1.7]{rydh_submersion_and_descent}.
\end{proof}

%

\begin{remark}
Examples of surjective morphisms that are universally submersive in the
constructible
topology are quasi-compact morphisms, universally open morphisms and
morphisms locally of finite
presentation~\cite[Prop.~1.6]{rydh_submersion_and_descent}.
Therefore, by Proposition~\pref{P:quotients-misc-props}, a universal Zariski
quotient is a topological quotient in the following cases:
\begin{enumerate}
\item $q$ is quasi-compact; 
\item $q$ is locally of finite presentation;
\item $q$ is universally open; 
\item $s$ is proper (then $q$ is universally closed with quasi-compact fibers,
hence quasi-compact); or
\item $s$ is universally open (then $q$ is universally open).
\end{enumerate}
In most applications, $s$ is universally open or even flat and locally of
finite presentation so that $q$ is also universally open.
\end{remark}

\begin{remark}
The definitions of topological and geometric quotients given above are not
standard but generalize other common definitions. In~\cite{kollar_quotients},
the conditions on topological and geometric quotients are slightly
stronger: all algebraic spaces are assumed to be locally noetherian, and
topological and geometric quotients are required to be locally of finite type.

In~\cite{GIT} the word topological is not defined explicitly but if we take
topological to mean the first three conditions of a geometric quotient in
\cite[Def.~0.6]{GIT}, taking into account that iii) should be universally
submersive, then a topological quotient is what here is called a universal
Zariski quotient.
\end{remark}

\begin{remark}
The condition, for a strongly topological quotient, that $j_{/Y}$ should be
universally submersive is natural. Indeed, this ensures that the equivalence
relation $X\times_Y X\inj X\times_S X$ has the quotient topology induced from
the groupoid. When $Y$ is a \emph{scheme}, or more generally a locally
separated algebraic space, the monomorphism $X\times_Y X\inj X\times_S X$ is an
immersion. In this case, the topology on $X\times_Y X$ is induced by
$X\times_S X$ and does not necessarily coincide with the quotient topology
induced by $R\surj X\times_Y X$.

If the groupoid has \emph{proper} diagonal $\map{j}{R}{X\times_S X}$, then every
topological quotient is strongly topological. This explains why proper group
actions are more tractable. An important class of strongly topological
quotients are topological quotients such that $\map{j_{/Y}}{R}{X\times_Y X}$ is
proper.
\end{remark}

\begin{remark}\label{R:geom=>cat}
A geometric quotient of schemes is always categorical in the category of
schemes, cf.~\cite[Prop.~0.1, p.~4]{GIT}, but not necessarily in the category
of algebraic spaces. As Koll\'ar mentions in~\cite[Rmk.~2.20]{kollar_quotients}
it is likely that every geometric quotient is categorical in the category of
\emph{locally separated} algebraic spaces. This is indeed the case, at least
for universally open quotients, as shown in
Theorem~\pref{T:(strong)-geom+univ-open-is-(strong)-cat}.

A natural condition, ensuring that a geometric quotient is categorical among
all algebraic spaces, is that the \emph{descent condition},
cf.\ Definition~\pref{D:descent-condition}, is fulfilled. Universally
open \emph{strongly geometric} quotients satisfy the descent condition
and are hence categorical,
cf.~Theorem~\pref{T:(strong)-geom+univ-open-is-(strong)-cat}.
\end{remark}

\begin{remark}\label{R:cat+top=>geom}
Conversely a (strongly) topological and uniformly categorical quotient is
(strongly) geometric. This is easily seen by considering equivariant maps $X\to
\A{1}_\Z$. Koll\'ar has also shown that if $G$ is an affine group, flat and
locally of finite type over $S$ acting \emph{properly} on $X$,
cf.\ Remark~\pref{R:proper-action}, such that a topological quotient exists,
then a geometric quotient exists~\cite[Thm.~3.13]{kollar_quotients}.
\end{remark}

\begin{proposition}\label{P:quotients-and-base-change}
Let $\groupoid{R}{X}$ be a groupoid and let $\map{q}{X}{Y}$ be an
equivariant morphism. Furthermore, let $\map{f}{Y'}{Y}$ be a morphism and let
$\map{q'}{X'}{Y'}$ be the pull-back of $q$ along $f$.
\begin{enumerate}
\item If $q$ is a topological quotient, then so is $q'$.
\item Assume that $f$ is flat. If $q$ is a geometric quotient, then so is $q'$.
\item\label{PI:qbc:fpqc-descent}
Assume that $f$ is covering in the fpqc topology (e.g., let $f$ be
faithfully flat and quasi-compact or locally of finite presentation).
If $q'$ is a topological (resp.\ geometric, resp.\ universal geometric)
quotient, then so is $q$.
\end{enumerate}
In particular, a geometric quotient is always uniform. The statements remain
valid if we replace ``topological'' and ``geometric'' with ``strongly
topological'' and ``strongly geometric''.
\end{proposition}
\begin{proof}
Topological and strongly topological quotients are universal by definition.
Part~\ref{PI:qbc:fpqc-descent}
for topological and strongly topological quotients follows immediately as a
morphism covering in the fpqc topology is submersive in both the Zariski and
the constructible
topology.
As exactness of a sequence of sheaves is preserved by flat base change and
can be checked fpqc locally, the statements on geometric quotients follow
by considering the sequence~\eqref{E:geometric-quot-seq}.
\end{proof}

\end{subsection}


\begin{subsection}{Separation properties}

Even if $X$ is separated a quotient $Y$ need not be. A sufficient criterion is
that the groupoid has proper diagonal and a precise condition, for
\emph{schemes}, is that the image of the diagonal is closed.

\begin{remark}\label{R:proper-action}
Consider the following properties of the diagonal of a group\-oid
$\groupoid{R}{X}$:
\begin{enumerate}
\item\label{RI:proper-action} the diagonal $\map{j}{R}{X\times_S X}$ is proper;
\item the diagonal $\map{j_{/Y}}{R}{X\times_Y X}$, with respect
to an equivariant morphism $\map{q}{X}{Y}$, is proper;
\item the diagonal $\map{j}{R}{X\times_S X}$ is quasi-compact;
\item the stabilizer $\map{\stab(X)}{j^{-1}\bigl(\Delta(X)\bigr)}{X}$ is
proper; and
\item the diagonal $\map{j}{R}{X\times_S X}$ is a monomorphism.
\end{enumerate}
If the groupoid is flat and locally of finite presentation, then, by the
cartesian diagrams in Remark~\pref{R:stacks-vs-groupoids:diagonal/inertia},
these properties correspond to the
following separation properties of the algebraic stack $\stX=[\groupoid{R}{X}]$:
\begin{enumerate}
\item $\stX$ is separated (over $S$);
\item $\stX\to Y$ is separated;
\item $\stX$ has quasi-compact diagonal (over $S$);
\item the inertia stack $I_\stX\to \stX$ is proper; and
\item $\stX$ is an algebraic space.
\end{enumerate}
If $\map{q}{X}{Y}$ is a topological quotient, then, by
Proposition~\pref{P:proper-actions} below, these properties imply that:
\begin{enumerate}
\item $Y$ is separated (over $S$); and
\addtocounter{enumi}{1}
\item $Y$ is quasi-separated (over $S$).
\end{enumerate}
If $(R,X)$ is the groupoid associated to a group action, then the group action
is called \emph{proper} if~\ref{RI:proper-action} holds. There is also the
notion of a
\emph{separated} group action which means that the diagonal $j$ has closed
image.
\end{remark}

\begin{proposition}\label{P:proper-actions}
Let $\groupoid{R}{X}$ be a groupoid and let $\map{q}{X}{Y}$ be a
topological quotient.
\begin{enumerate}
\item $j$ is proper if and only if $Y$ is separated and $j_{/Y}$ is
proper.
\item $j$ is quasi-compact if and only if $Y$ is quasi-separated
and $j_{/Y}$ is quasi-compact.
\item If $j_{/Y}$ is proper, then the stabilizer
is proper.\label{I:prop-stab}
\item If $Y$ is locally separated, then $Y$ is separated if and only if
the image of $j$ is closed.
\end{enumerate}
\end{proposition}
\begin{proof}
As $q$ is a topological quotient we have that $q$ and $q^\cons$ are universally
submersive and that $j_{/Y}$ is surjective. The statements then easily follow
from~\cite[Prop.~1.7]{rydh_submersion_and_descent} and the cartesian diagram
\begin{equation*}\raisebox{1.5 cm}{\xymatrix{%
{j\colon R}\ar[r]^-{j_{/Y}} &
X\times_Y X\ar@{(->}[r]\ar[d] & X\times_S X\ar[d]^{q\times q} \\
& Y\ar@{(->}[r]^-{\Delta_{Y/S}} & Y\times_S Y.\ar@{}[ul]|\square}}\qedhere
\end{equation*}
\end{proof}

\begin{remark}[{\cite[Cor.~5.2]{conrad_coarsespace}}]
\label{R:proper-stabilizer}
Assume that there exists an equivariant morphism $\map{q}{X}{Y}$ with respect
to the groupoid $(R,X)$
such that the diagonal $j_{/Y}$ is proper. Then the groupoid has a proper
stabilizer. Moreover, if $\map{r}{X}{Z}$ is a
\emph{categorical} quotient, then it follows that $j_{/Z}$ is proper as well. In
particular, if $\stX$ is an algebraic stack admitting a separated morphism
$\map{r}{\stX}{Y}$ to an algebraic space $Y$, then $\stX$ has a proper inertia
stack and the categorical quotient $\map{q}{\stX}{Z}$, if it exists, is
separated.

The Keel--Mori theorem asserts that,
conversely, if $\groupoid{R}{X}$ is a flat groupoid locally of finite
presentation such that the stabilizer map $j^{-1}\bigl(\Delta(X)\bigr)\to X$
is \emph{finite}, then there exists a geometric and categorical
quotient $X\to Z$ and $j_{/Z}$ is proper, cf.\ Theorem~\pref{T:KM-general}.
\end{remark}

As we will be particularly interested in finite groupoids, we make the following
observations.

\begin{proposition}\label{P:proper-actions-for-proper-groups}
Let $\groupoid{R}{X}$ be a \emph{proper} groupoid, i.e., $s$ and $t$ are
proper, and let $\map{q}{X}{Y}$ be a topological quotient.
\begin{enumerate}
\item\label{PI:proper-act:j} The diagonal $j$ is proper (resp.\ quasi-compact)
if and only if $X$ is
separated (resp.\ quasi-separated).
\item\label{PI:proper-act:j_Y} The diagonal $j_{/Y}$ is proper
(resp.\ quasi-compact) if and only if $q$
is separated (resp.\ quasi-separated).
\end{enumerate}
\end{proposition}
\begin{proof}
As $s$ is separated the section $\map{e}{X}{R}$ is closed. Thus
$\map{\Delta_{X/S}=j\circ e}{X\inj R}{X\times_S X}$ is proper
(resp.\ quasi-compact) if $j$ is so. Conversely, if $X$ is separated
(resp.\ quasi-separated), then, as $s=\pi_1\circ j$ is proper,
it follows that $j$ is proper (resp.\ quasi-compact).
\ref{PI:proper-act:j_Y} follows from \ref{PI:proper-act:j} considering the
$S$-groupoid as a $Y$-groupoid.
\end{proof}





\end{subsection}


\begin{subsection}{Free actions}

\begin{definition}
We say that a groupoid $\groupoid{R}{X}$ is an \emph{equivalence
relation} if $\map{j}{R}{X\times_S X}$ is a monomorphism. Let $X$ be an
algebraic space over $S$ with an action of a group scheme $G\to S$. We say that
$G$ acts
\emph{freely} if the associated groupoid is an equivalence relation, i.e., if
the morphism $\map{j}{G\times_S X}{X\times_S X}$ is a monomorphism.
\end{definition}

\begin{theorem}\label{T:equivalence-relation}
  Let $\groupoid{R}{X}$ be a flat and locally finitely presented
  equivalence relation. Then there is a universal
  strongly geometric and categorical quotient $\map{q}{X}{X/R}$ in the category
  of algebraic spaces. Furthermore, $q$ is the quotient of the equivalence
  relation in the category of fppf sheaves. Hence, $q$ is flat and locally of
  finite presentation.
\end{theorem}
\begin{proof}
Let $X/R$ be the quotient sheaf of the equivalence relation
$\groupoid{R}{X}$ in the fppf topology. Then $X/R$ is an algebraic space
by~\cite[Cor.~6.3]{artin_versal_def_alg_stacks} as explained
in~\cite[Cor.~10.4]{laumon} and~\cite[Thm.~04S6]{stacks-project}.
As $X/R$ is a categorical quotient in the category
of fppf sheaves, it is a categorical quotient in the category
of algebraic spaces. As taking the quotient sheaf commutes with arbitrary base
change, it is also a universal categorical quotient.

The quotient $\map{q}{X}{X/R}$ is flat and locally of finite presentation and
thus universally submersive in both
the Zariski and the constructible topology. As
$\map{j_{/(X/R)}}{R}{X\times_{X/R} X}$ is an isomorphism, we thus have that $q$
is a
strongly topological quotient. It is then a universal strongly geometric
quotient by Remark~\pref{R:cat+top=>geom}.
\end{proof}

\begin{corollary}
  Let $G\to S$ be a group scheme, flat and locally of finite presentation over
  $S$. Let
  $X\to S$ be an algebraic space with a free action of $G$. Then there is a
  universal strongly geometric and categorical quotient $\map{q}{X}{X/G}$ in
  the category of algebraic spaces and $q$ is flat and locally of finite
  presentation.
\end{corollary}
\begin{proof}
This follows immediately from Theorem~\pref{T:equivalence-relation}.
\end{proof}

\end{subsection}

\end{section}


\begin{section}{Fixed-point reflecting morphisms and the descent condition}
\label{S:fpr-and-descent}

Let $X$ and $Y$ be algebraic spaces with actions of a group $G$ such that
geometric quotients $X/G$ and $Y/G$ exist. If $\map{f}{X}{Y}$ is any \etale{}
$G$-equivariant morphism, then, in general, the induced morphism
$\map{f/G}{X/G}{Y/G}$ is not \etale{}. The notion of fixed-point reflecting
morphisms was introduced to remedy this problem. Under a mild assumption, $f/G$
is
\etale{} if $f$ is fixed-point reflecting. Moreover,
the existence of $X/G$ then follows from the existence of $Y/G$ and
$X=X/G\times_{Y/G} Y$, that is, $f$ is \emph{strongly \etale{}}.
The mild assumption that we need is that $\map{q}{Y}{Y/G}$ satisfies effective
descent for \etale{} morphisms. This is almost always the case, cf.\ 
Appendix~\ref{A:descent}.

An interpretation of this section in terms of stacks is given
in~\S\ref{S:cms-stacks}.

\begin{remark}
Recall that the stabilizer of a point
$\map{x}{\Spec(k(x))}{X}$ is the fiber
$\stab(x)=j^{-1}(x,x)=\stab(X)\times_X \Spec(k(x))$, which is
a group scheme over $k(x)$, cf.\ Notation~\pref{N:stabilizer}.
Assume that $\groupoid{R}{X}$ is the groupoid
associated to the action of a group $G$ on $X$. Then
$\stab(x)\subseteq G$ is the subgroup of elements $g\in G$ such that $g(x)=x$
(as morphisms $\Spec(k(x))\to X$, not as points in $|X|$).
\end{remark}

\begin{definition}\label{D:square}
Let $\map{f}{(R_X,X)}{(R_Y,Y)}$ be a morphism of groupoids. We say that $f$
is \emph{square} if the two commutative diagrams
$$\vcenter{
\xymatrix{
R_X\ar[r]^-f\ar[d]^{s} & R_Y\ar[d]^{s}\\
X\ar[r]^-f & Y}}
\quad\text{and}\quad
\vcenter{\xymatrix{
R_X\ar[r]^-f\ar[d]^{t} & R_Y\ar[d]^{t}\\
X\ar[r]^-f & Y}}
$$
are cartesian.
\end{definition}

Note that if one of the diagrams above is cartesian, then so is the other
diagram, since the involution of the groupoid exchanges the two diagrams.

\begin{definition}\label{D:fpr}
Let $\map{f}{(R_X,X)}{(R_Y,Y)}$ be a morphism of groupoids. We say that $f$ is
\emph{fixed-point reflecting} at $x\in |X|$, abbreviated fpr, if
the canonical
morphism of stabilizers $\stab(x)\to \stab(f\circ x)$
is an
isomorphism for some, or equivalently every, representative
$\map{x}{\Spec(k)}{X}$.
We let $\fpr(f)\subseteq |X|$ denote the subset over which
$f$ is fixed-point reflecting.
\end{definition}

The fixed-point reflecting set $\fpr(f)\subseteq |X|$ is
$R_X$-stable. Indeed, if $\map{r}{\Spec(k)}{R_X}$ is a point, then there
is a commutative diagram
$$\xymatrix{%
\stab(s\circ r)\ar[d]^{\iso}\ar[r] & \stab(s\circ f\circ r)\ar[d]^{\iso}\\
\stab(t\circ r)\ar[r] & \stab(t\circ f\circ r)}$$
so that $f$ is fpr at $s\circ r$ if and only if $f$ is fpr at $t\circ r$.
Some authors prefer the terminology \emph{stabilizer preserving} rather than
fixed-point reflecting.

\begin{remark}\label{R:fpr-from-base-change}
If $X\to Z$ is an equivariant morphism with respect to the groupoid $(R_X,X)$
and $Z'\to Z$ is any
morphism, then $(R_X\times_Z Z',X\times_Z Z') \to (R_X,X)$ is fpr.
\end{remark}

The following proposition sheds some light over the importance of proper
stabilizer.

\begin{proposition}\label{P:fpr-of-unramified-is-open}
Let $\map{f}{(R_X,X)}{(R_Y,Y)}$ be a square morphism of groupoids such that
$\map{f}{X}{Y}$ is \emph{unramified}. If the stabilizer $\stab(Y)\to Y$ is
\emph{universally closed}, then the subset $\fpr(f)$ is open in $X$.
\end{proposition}
\begin{proof}
There are cartesian diagrams
$$\xymatrix{
\stab(X)\ar@{(->}[r]\ar[d] & \stab(Y)\times_Y X\ar@{(->}[r]\ar[d]
   & R_X\ar[r]\ar[d]^{j} & R_Y\ar[dd]^{j} \\
X\ar@{(->}[r]^-{\Delta_{X/Y}} & X\times_Y X\ar@{(->}[r]\ar^{\pi_2}[d]\ar@{}[ul]|\square
   & X\times_S X\ar[d]\ar@{}[ul]|\square \\
& X\ar@{(->}[r] & Y\times_S X\ar[r]\ar@{}[ul]|\square
  & Y\times_S Y.\ar@{}[uul]|\square
}$$
The morphism $\map{\varphi}{\stab(Y)\times_Y X}{X}$, which is the second column
above, is closed.
A point $x\in |X|$ is fpr if and only if $x$ is not in $\varphi\bigl({|\stab(Y)\times_Y X|}\setminus |\stab(X)|\bigr)$.
As $f$ is unramified, $\Delta_{X/Y}$ is an open
immersion and hence $|\stab(Y)\times_Y X|\setminus |\stab(X)|$ is closed. Thus
$\fpr(f)$ is
open.
\end{proof}


\begin{definition}\label{D:descent-condition}
Let $\groupoid{R_X}{X}$ be a groupoid and let $\map{q}{X}{Z_X}$ be a
topological quotient. We say that the quotient $q$ satisfies the
\emph{descent condition} if
for any \etale{}, square and fixed-point reflecting morphism of groupoids
$\map{f}{(R_W,W)}{(R_X,X)}$, there exists an algebraic space $Z_W$ and a
cartesian square
\begin{equation}\label{E:descent-condition}
\vcenter{\xymatrix{
W\ar[r]^f\ar[d] & X\ar[d]^q \\
Z_W\ar[r] & Z_X\ar@{}[ul]|\square
}}\end{equation}
where $Z_W\to Z_X$ is \emph{\etale}. We say that $q$ satisfies the
\emph{weak descent condition} if the descent condition holds when restricted
to morphisms $f$ such that there is a cartesian square
$$\xymatrix{
W\ar[r]^f\ar[d] & X\ar[d]^q \\
Q'\ar[r] & Q\ar@{}[ul]|\square
}$$
where $X\to Q$ is an equivariant morphism to a \emph{locally separated}
algebraic
space $Q$ and $Q'\to Q$ is an \etale{} morphism.

Moreover, we say that $q$ satisfies the (weak) descent condition for $\Et_P(X)$
for $P\in\{\qc,\qsep,\cons,\sep,\sep+\qc,\finite\}$ if the (weak) descent
condition holds when restricted to \etale{} morphisms $W\to X$ with property
$P$,
cf.\ Appendix~\ref{A:descent}.
\end{definition}

\begin{remark}
If a space $Z_W$ exists as above, then $W\to Z_W$ is a topological quotient.
If $q$ is a geometric quotient, then $W\to Z_W$ is also a geometric quotient by
Proposition~\pref{P:quotients-and-base-change}.
A morphism $\map{f}{W}{X}$ is \emph{strongly \etale{}} if the induced
morphism on categorical quotients $Z_W\to Z_X$ is \etale{} and the
diagram~\eqref{E:descent-condition} is cartesian~\cite[p.~198]{GIT}. Thus, the
descent condition ensures that every fixed-point reflecting \etale{} square
morphism is strongly \etale{}.
In~\cite[2.4]{keel_mori_quotients} the existence of $Z_W$ is not a part of the
descent condition.
\end{remark}

For most purposes, such as in the following proposition, it is sufficient to
restrict the descent condition to \emph{separated} \etale{} morphisms.

\begin{proposition}[{\cite[Cor.~2.15]{kollar_quotients}}]\label{P:geom+dc=>cat}
Let $\groupoid{R}{X}$ be a groupoid and let $\map{q}{X}{Z}$ be a
geometric quotient satisfying the descent condition (resp.\ the weak descent
condition) for the category $\Et_{\sep}(X)$ of separated \etale{} morphisms
$W\to X$.
Then $q$ is a categorical quotient (resp.\ a categorical quotient
among locally separated algebraic spaces).

Moreover, it is enough that $q$ satisfies the (weak) descent condition for
$\Et_{\sep+\qc}(X)$ if we restrict the discussion to the category of
quasi-separated algebraic spaces.
\end{proposition}
\begin{proof}
Let $Q$ be an algebraic space (resp.\ a locally separated algebraic space)
and let $\map{r}{X}{Q}$ be an equivariant morphism.
We have to prove that there is a unique morphism $\map{g}{Z}{Q}$
such that $r=g\circ q$. As geometric quotients commute with open immersions, we
can assume that $Q$ is quasi-compact. Let $Q'\to Q$ be an \etale{} presentation
with $Q'$ an affine scheme. Let $X'=X\times_Q Q'$ so that $X'\to X$ is separated
and \etale{} and also quasi-compact if $Q$ is quasi-separated. As $q$ satisfies
the descent
condition (resp.\ the weak descent condition), there is a geometric quotient
$\map{q'}{X'}{Z'}$ such that $X'=X\times_Z Z'$.

As $Q'$ is affine, the morphism
$X'\to Q'$ is determined by $\Gamma(Q')\to\Gamma(X')$. Moreover, as $q'$ is
geometric, we have that the image of $\Gamma(Q')$ lies in
$\Gamma(Z')=\Gamma(X')^{R'}$. The induced homomorphism
$\Gamma(Q')\to\Gamma(Z')$ gives a morphism $\map{g'}{Z'}{Q'}$ such that
$r'=g'\circ q'$ and this is the only morphism $g'$ with this property. By
\etale{} descent, the morphism $g'$ descends to a unique morphism
$\map{g}{Z}{Q}$ such that $r=g\circ q$.
\end{proof}

\begin{definition}
Let $\map{f}{X}{Y}$ be a morphism of algebraic spaces. Let $W\subseteq |X|$
be a subset. We say that $f$ is universally injective (resp.\ universally
bijective, resp.\ a universal homeomorphism, resp.\ universally submersive)
over $W$ if for any cartesian diagram
$$\xymatrix{%
X'\ar[r]^{g'}\ar[d]_{f'} & X\ar[d]^f \\
Y'\ar[r]^g & Y\ar@{}[ul]|\square}$$
we have that the map $\map{f'|_{g'^{-1}(W)}}{g'^{-1}(W)}{|Y'|}$ is injective
(resp.\ bijective, resp.\ a homeomorphism, resp.\ submersive).
\end{definition}

\begin{lemma}\label{L:univ-homeo-over-subset-of-unramified}
Let $\map{f}{X}{Y}$ be an unramified morphism of algebraic spaces. Let
$W\subseteq |X|$ be a subset such that $f$ is a universal homeomorphism over
$W$. Then $W\subseteq |X|$ is open.
\end{lemma}
\begin{proof}
Let $\map{\pi_1,\pi_2}{X\times_Y X}{X}$ be the projections and let
$W_1=\pi_1^{-1}(W)$. Then, by assumption, we have that
$\map{\pi_2|_{W_1}}{W_1}{|X|}$ is a
homeomorphism. The diagonal $\map{\Delta_f}{X}{X\times_Y X}$ is an open
immersion. Thus $W_1\cap \Delta_f(X)$ is open in $W_1$. It follows that
$W=\pi_2\bigl(W_1\cap \Delta_f(X)\bigr)$ is open in $X$.
\end{proof}

We will now describe the descent condition in terms of effective descent of
\etale{} morphisms.

\begin{setup}\label{SE:fpr-gamma}
Let $\map{f}{(R_W,W)}{(R_X,X)}$ be an \etale{} and square morphism of groupoids
and let $X\to Z_X$ be a topological quotient. Consider the following diagram
\begin{equation}\label{E:phi-diagram}
\vcenter{\xygraph{ !{0;<35mm,0cm>:<0cm,1cm>:: \labelmargin+{1mm}}
[d] *+{R_W}="R"
[ur] *+{W\times_{Z_X} X}="WX" 
[dd] *+{X\times_{Z_X} W}="XW"
[ur] *+{X\times_{Z_X} X.}="XX"
"R" : "WX" ^-{(s,f\circ t)}
"R" : "XW" _-{(f\circ s,t)}
"WX" :@{-->} "XW" ^{\varphi}
"WX" : "XX" ^-{f\times \id{X}}
"XW" : "XX" _-{\id{X}\times f}
}}
\end{equation}
As $X\to Z_X$ is topological, the two morphisms on the left are surjective.
To give a $X\times_{Z_X} X$-morphism
$$\map{\varphi}{W\times_{Z_X} X}{X\times_{Z_X} W}$$
is equivalent to specifying its graph
$$\Gamma_\varphi\inj W\times_{Z_X} W\iso
(W\times_{Z_X} X)\times_{X\times_{Z_X} X}(X\times_{Z_X} W),$$
which is an open subspace as $\map{f}{W}{X}$ is \etale{}.
The morphism $\varphi$ makes the diagram above commutative if and only if
$$\map{j_{W/Z_X}=(s,t)}{R_W}{W\times_{Z_X} W}$$
factors through $\Gamma_\varphi$. Let $\Gamma_{W/Z_X}\subseteq
|W\times_{Z_X} W|$
denote the image of $j_{W/Z_X}$. Then $\varphi$ makes the diagram
commutative if and only if $\Gamma_{W/Z_X}\subseteq \Gamma_\varphi$, or
equivalently, if and only if $\Gamma_{W/Z_X}=\Gamma_\varphi$ since
$(s,f\circ t)$
is surjective. In particular, there is at most one morphism $\varphi$ making
the diagram commutative.
\end{setup}

\begin{proposition}\label{P:descent-data}
Let $\map{f}{W}{X}$ and $X\to Z_X$ be as in Setup~\pref{SE:fpr-gamma}.
There is a one-to-one correspondence between, on the one hand, topological
quotients
$W\to Z_W$
together with an \etale{} morphism $Z_W\to Z_X$ such that the diagram
\eqref{E:descent-condition}
is cartesian and, on the other hand, \emph{effective} descent data $\varphi$ for
$\map{f}{W}{X}$ fitting into the diagram~\eqref{E:phi-diagram}.
In particular, as there is at most one such morphism $\varphi$, there is at most
one quotient $Z_W$ in the descent condition~\pref{D:descent-condition}.

Assume that $f$ is fixed-point reflecting. Then the following holds.
\begin{enumerate}
\item \label{P:dd:fpr-univ-bijective}
The morphisms
$$\map{f\times \id{W}}{W\times_{Z_X} W}{X\times_{Z_X} W}$$
$$\map{\id{W}\times f}{W\times_{Z_X} W}{W\times_{Z_X} X}$$
are universally bijective over the subset $\Gamma_{W/Z_X}$.
\item \label{P:dd:univ-sub}
The subset $\Gamma_{W/Z_X}$ is open if and only if $f\times \id{W}$ (or
$\id{W}\times f$) is universally submersive over $\Gamma_{W/Z_X}$.
\item \label{P:dd:open}
An isomorphism $\varphi$ fitting into the diagram~\eqref{E:phi-diagram} exists
if and only if $\Gamma_{W/Z_X}$ is open and then $\Gamma_{W/Z_X}$ is the graph
of $\varphi$.
\end{enumerate}
\end{proposition}
\begin{proof}
There is a one-to-one correspondence between \etale{} morphisms $Z_W\to Z_X$ and
\etale{} morphisms $\map{f}{W}{X}$ together with an effective descent datum
$\varphi$,
i.e., an isomorphism $\map{\varphi}{W\times_{Z_X} X} {X\times_{Z_X} W}$ over
$X\times_{Z_X} X$ satisfying the cocycle condition. This is because $q$ is
universally submersive and hence a
morphism of descent for \etale{} morphisms, cf.\ Proposition~\pref{P:descent}.
Given a quotient $Z_W$ of $\groupoid{R_W}{W}$ together with an \etale{} morphism
$Z_W\to Z_X$ as in
Definition~\pref{D:descent-condition}, we have that the corresponding
isomorphism $\varphi$ is the composition of the canonical isomorphisms
$W\times_{Z_X} X\iso W\times_{Z_W} W\iso X\times_{Z_X} W$ and that $\varphi$
fits into the commutative diagram~\eqref{E:phi-diagram}.

\ref{P:dd:fpr-univ-bijective}:
By symmetry, it is enough to show that the first morphism is universally
bijective over $\Gamma_{W/Z_X}$. Let $\map{(x_1,w_2)}{\Spec(k)}{X\times_{Z_X}
  W}$ be a point. We have to show that there is exactly one lifting to
$W\times_{Z_X} W$ in the image of $R_W$. As $R_W\to X\times_{Z_X} W$ is
surjective, we can enlarge $k$ and assume that there exists a lifting
$\map{r}{\Spec(k)}{R_W}$.
We are free to enlarge $k$, so it is enough to show that any two
liftings to $R_W$ have the same image in $W\times_{Z_X} W$.

Let $(w_1,w_2)$ be the image of $r$
in $W\times_{Z_X} W$. The liftings of $(w_1,w_2)$ to $R_W$ are in bijection,
via conjugation by $r$, to $k$-points of the stabilizer group scheme
$\stab(w_1)$.
On the other hand, as $\map{f}{W}{X}$ is square, we have that liftings of
$(x_1,w_2)$ to $R_W$ are in bijection with $k$-points of $\stab(x_1)$. As $f$ is
fixed-point reflecting, these two sets are canonically bijective so that every
lifting of $(x_1,w_2)$ to $R_W$ comes from a lifting of $(w_1,w_2)$ to $R_W$.
Hence $(w_1,w_2)$ is the unique lifting of $(x_1,w_2)$ in the image of $R_W$.

\ref{P:dd:univ-sub} follows from \ref{P:dd:fpr-univ-bijective} and
Lemma~\pref{L:univ-homeo-over-subset-of-unramified}.

\ref{P:dd:open}: In the setup we have seen that if a morphism $\varphi$ exists,
then necessarily $\Gamma_{W/Z_X}=\Gamma_\varphi$ is open. Conversely, if
$\Gamma_{W/Z_X}$ is open, then $\Gamma_{W/Z_X}\to W\times_{Z_X} X$ and
$\Gamma_{W/Z_X}\to X\times_{Z_X} W$ are universally bijective
by~\ref{P:dd:fpr-univ-bijective} and hence
isomorphisms~\cite[Exp.~IX, Cor.~1.6]{sga1}. It follows that $\Gamma_{W/Z_X}$ is
the graph of an isomorphism $\varphi$ fitting into
diagram~\eqref{E:phi-diagram}.
\end{proof}

\begin{corollary}\label{C:dc-and-base-change}
Let $\groupoid{R}{X}$ be a groupoid and let $\map{q}{X}{Z}$ be a
topological quotient. Also let $\map{g}{Z'}{Z}$ be a morphism and let
$\map{q'}{X'}{Z'}$ be the pull-back of $q$ along $g$. Then $q'$ is a
topological quotient by Proposition~\pref{P:quotients-and-base-change}. Let
$P\in \{\qc,\qsep,\cons,\sep,\sep+\qc,\finite\}$ be a property of \etale{}
morphisms.
\begin{enumerate}
\item\label{CI:dc-bc:stab} Assume that $g$ is \etale{} and has property $P$.
If $q$ satisfies the descent
condition with respect to $\Et_P$, then so does $q'$.
\item\label{CI:dc-bc:descent} Assume that $g$ is covering in the fppf topology.
Then $q$
satisfies the descent condition for $\Et_P$, if $q'$
satisfies the descent condition for $\Et_P$.
\end{enumerate}
\end{corollary}
\begin{proof}
To prove \ref{CI:dc-bc:stab}, let $\map{f'}{W'}{X'}$ be an \etale{}, square and
fpr morphism of groupoids with property $P$. Then $W'\to X'\to X$ is also
\etale{}, square and fpr with property $P$.
As $q$ satisfies the descent condition there is a
topological quotient $W'\to Z_{W'}$ and an \etale{} morphism $Z_{W'}\to Z$.

We need to construct a $Z$-morphism $Z_{W'}\to Z'$ such that
$W'=Z_{W'}\times_{Z'} X'$, that is, we need to descend the morphism $f'$ along
$q$. The local description of this morphism is $\map{q'\circ f'}{W'}{Z'}$ and
since $q$ is a morphism of descent, it is enough to verify that
$\map{q'\circ f'\circ\pi_1,q'\circ f'\circ\pi_2}
{W'\times_{Z_{W'}} W'}{Z'}$ coincide.
Equivalently, it is enough to show
that $\map{f'\times f'}{W'\times_{Z_{W'}} W'}{X'\times_Z X'}$ factors through
the open subspace $X'\times_{Z'} X'$. This follows from the commutative diagram
$$\xymatrix{%
R_{W'}\ar[r]\ar[d] & W'\times_{Z_{W'}} W'\ar[r]^-{f'\times f'} & X'\times_Z X'\\
R_{X'}\ar[r] & X'\times_{Z'} X'\ar@{(->}[ur]}$$
%
%
as $R_{W'}\to W'\times_{Z_{W'}} W'$ is surjective.

\ref{CI:dc-bc:descent} follows from an easy application of fppf descent, taking
into account that,
by Proposition~\pref{P:descent-data}, the quotient $Z_W$ figuring in the
descent condition is unique.
\end{proof}

The following theorem generalizes~\cite[Thm.~2.14]{kollar_quotients}.
In Koll\'{a}r's result, the group action has to be proper, i.e., the diagonal
$\map{j}{R}{X\times_S X}$ has to be proper, so that every topological quotient
is strongly topological.

\begin{theorem}\label{T:descent-condition-holds}
Let $\groupoid{R}{X}$ be a groupoid and let $\map{q}{X}{Z}$ be a
topological
quotient (resp.\ a strongly topological quotient) such that $q$ satisfies
\emph{effective} descent for $\Et_P(X)$. Then $q$ satisfies the weak
descent condition (resp.\ the descent condition) for $\Et_P(X)$.
\end{theorem}
\begin{proof}
Let $\map{f}{W}{X}$ be an \etale{}, square and fpr morphism with property~$P$.
As $q$
satisfies effective descent for $f$ it is, by Proposition~\pref{P:descent-data},
enough to show that the image $\Gamma_{W/Z}$ of
$\map{j_{W/Z}}{R_W}{W\times_{Z} W}$ is open, or equivalently, that
$\map{f\times \id{W}}{W\times_{Z} W}{X\times_{Z} W}$
is universally submersive over $\Gamma_{W/Z}$. If $q$ is a strongly topological
quotient, then
$(f\times \id{W})\circ j_{W/Z}$ is universally submersive since it is the
pull-back of the universally submersive morphism
$\map{j_{/Z}}{R_X}{X\times_Z X}$.
Thus $f\times \id{W}$ is universally submersive over $\Gamma_{W/Z}$. This shows
that $q$ satisfies the descent condition.

Now, let $\map{q}{X}{Z}$ be a topological quotient that satisfies effective
descent for $\Et_P(X)$ and let $\map{r}{X}{Q}$ be an
equivariant morphism such that $Q$ is a locally separated algebraic space. Then $X\times_Q
X\inj X\times_S X$ is an \emph{immersion}. Thus
$$\map{\pi_1}{(X\times_Z X)\times_{X\times_S X} (X\times_Q X)}{X\times_Z X}$$
is an immersion. Moreover, as $R_X\to X\times_Z X$ is surjective we have that
$\pi_1$ is surjective. Hence we obtain a monomorphism $(X\times_Z X)_\red\to
X\times_Q X$ over $X\times_S X$.

Let $Q'\to Q$ be an \etale{} morphism with property $P$ and let
$W=X\times_Q Q'$ and
$R_W=R_X\times_Q Q'$ so that $(R_W,W)\to (R_X,X)$ is fixed-point reflecting and
square \etale{}, cf.\ Remark~\pref{R:fpr-from-base-change}.
To show that $q$ satisfies the weak descent condition,
it is by Proposition~\pref{P:descent-data} enough to show that
$\Gamma_{W/Z}\subseteq |W\times_Z W|$ is open. We have a cartesian diagram
$$\xymatrix{%
(R_W)_\red\ar[dd]\ar[r]_-{(s,f\circ t)}\ar@/^1pc/[rr]^{(s,t)}
 & (W\times_Z X)_\red\ar@{(->}[r]\ar@{(->}[d]^{\psi}
 & W\times_{Q'} W\ar@{(->}[d]\ar[r] & Q'\ar@{(->}[d]^{\Delta_{Q'/Q}}\\
 & (W\times_Z W)_\red\ar@{(->}[r]\ar[d]
 & W\times_Q W\ar[d]\ar[r]\ar@{}[ul]|\square
 & Q'\times_{Q} Q'\ar@{}[ul]|\square\\
(R_X)_\red\ar[r]^-{j_{/Z}}
 & (X\times_Z X)_\red\ar@{(->}[r]\ar@{}[uul]|{\square\quad\quad}
 & X\times_Q X.\ar@{}[ul]|\square}$$
As $\Delta_{Q'/Q}$ is an open immersion so is $\psi$ and as
$(s,f\circ t)$ is surjective it follows that the image of
$(j_{W/Z})_\red=\psi\circ (s,f\circ t)$ is the open subspace
$\Gamma_{W/Z}=\psi(W\times_Z X)$.
\end{proof}

\begin{theorem}\label{T:descent-condition-holds-for-univ-open/int}
Let $\groupoid{R}{X}$ be a groupoid and let $\map{q}{X}{Z}$ be a
topological (resp.\ strongly topological) quotient.
\begin{enumerate}
\item If $q$ is universally open, then $q$ satisfies the weak descent (resp.\
descent) condition for $\Et_{\qsep}(X)$.
\item If $q$ is proper or integral, then $q$ satisfies the weak descent
(resp.\ descent)
condition for $\Et(X)$.
\end{enumerate}
\end{theorem}
\begin{proof}
This follows from Theorem~\pref{T:descent-condition-holds} as $q$ is a morphism
of effective descent by Theorem~\pref{T:effective-descent}.
\end{proof}

The following theorem generalizes~\cite[Cor.~2.15]{kollar_quotients} and
answers the conjecture in~\cite[Rmk.~2.20]{kollar_quotients}.

\begin{theorem}\label{T:(strong)-geom+univ-open-is-(strong)-cat}
Let $\groupoid{R}{X}$ be a groupoid and let $\map{q}{X}{Z}$ be a
geometric (resp.\ strongly geometric) quotient. Consider the following
conditions on the quotient~$q$:
\begin{enumerate}
\myitem{1a}\label{TI:catq:uo} $q$ is universally open;
\myitem{1b}\label{TI:catq:p/i} $q$ is proper or integral;
\myitem{2a}\label{TI:catq:qc+ln} $q$ is quasi-compact and $Z$ is locally
noetherian; and
\myitem{2b}\label{TI:catq:qc+uc} $q$ is quasi-compact and universally
subtrusive (e.g., univ.\ closed).
\end{enumerate}
If~\ref{TI:catq:uo} or~\ref{TI:catq:p/i} holds, then $q$ is a categorical
quotient among locally
separated algebraic spaces
(resp.\ a categorical quotient).
If~\ref{TI:catq:qc+ln} or~\ref{TI:catq:qc+uc} holds, then $q$ is a categorical
quotient among locally
separated and quasi-separated algebraic spaces
(resp.\ a categorical quotient among quasi-separated algebraic spaces).
\end{theorem}
\begin{proof}
If~\ref{TI:catq:uo} or~\ref{TI:catq:p/i} holds, then $q$ is a morphism of
effective descent for \etale{}
and separated morphisms, by Theorem~\pref{T:effective-descent}.
Likewise, if~\ref{TI:catq:qc+ln} or~\ref{TI:catq:qc+uc} holds,
then $q$ is a morphism of
effective descent for \etale{}, quasi-compact and separated morphisms.
Therefore $q$ satisfies the weak descent
condition (resp.\ descent condition) for $\Et_{\sep}(X)$ or $\Et_{\sep+\qc}(X)$,
by Theorem~\pref{T:descent-condition-holds}, and is a categorical quotient by
Proposition~\pref{P:geom+dc=>cat}.
\end{proof}

\begin{definition}
An equivariant morphism is called a \GC{} quotient if it is a strongly geometric
quotient that satisfies
the descent condition for separated \etale{} morphisms uniformly. As a \GC{}
quotient is categorical by
Proposition~\pref{P:geom+dc=>cat}, we will speak about \emph{the} \GC{} quotient
when it exists.
\end{definition}

\begin{remark}\label{R:KM-quotients-dc}
  The definition of \GC{} quotient given by Keel and Mori differs slightly from
  ours. In~\cite{keel_mori_quotients} it simply means a geometric and uniform
  categorical quotient. However, every quotient $\map{q}{X}{Y}$ appearing
  in~\cite{keel_mori_quotients} is such that $\map{j_{/Y}}{R}{X\times_Y X}$ is
  proper, hence every quotient is strongly topological. As every
  groupoid in \loccit\ is flat and locally of finite presentation, it follows
  that every quotient is universally open. Every
  \GC{} quotient in~\cite{keel_mori_quotients} thus satisfies the descent
  condition for $\Et_\sep$ uniformly by
  Theorem~\pref{T:descent-condition-holds-for-univ-open/int}.
\end{remark}

\begin{theorem}[{\cite[Cor.~2.17]{kollar_quotients}}]
\label{T:existence-of-GC-by-fpr}
Let $\map{f}{(R_W,W)}{(R_X,X)}$ be a surjective, \etale{}, separated, square
and fpr
morphism of groupoids and let $Q=W\times_X W$. Assume that a \GC{} quotient
$W\to Z_W$ of $(R_W,W)$ exists. Then \GC{} quotients $Q\to Z_Q$ and $X\to Z_X$
exist and $Z_X$ is the quotient of the \etale{} equivalence relation
$\xymatrix@M=1pt{Z_Q\ar@<.5ex>[r] \ar@<-.5ex>[r] & Z_W}$. Moreover, the
natural squares of the diagram
\begin{equation}\label{E:existence-GC-fpr}
\vcenter{
\xymatrix{Q \ar@<.5ex>[r] \ar@<-.5ex>[r]\ar[d] & W \ar[r]^f\ar[d] & X\ar[d] \\
Z_Q \ar@<.5ex>[r] \ar@<-.5ex>[r] & Z_W \ar[r] & Z_X}}
\end{equation}
are cartesian.
\end{theorem}
\begin{proof}
Since $\map{f}{W}{X}$ is square, \etale{} and fixed-point reflecting so are the
two projections $\map{\pi_1,\pi_2}{Q=W\times_X W}{W}$. As $W\to Z_W$ satisfies
the descent condition, there exists topological quotients $Q\to
\left(Z_Q\right)_1$ and $Q\to \left(Z_Q\right)_2$ induced by the two
projections $\pi_1$ and $\pi_2$. Moreover, by
Corollary~\pref{C:dc-and-base-change} the quotients $Q\to
\left(Z_Q\right)_1$ and $Q\to\left(Z_Q\right)_2$ satisfy the descent condition
and it follows by Proposition~\pref{P:geom+dc=>cat} that
$\left(Z_Q\right)_1\iso\left(Z_Q\right)_2$ is the unique \GC{} quotient. The
two canonical morphisms $Z_Q\to Z_W$ are \etale{} and the corresponding squares
are cartesian.

We will now show that $\groupoid{Z_Q}{Z_W}$ is an equivalence relation. First
note that we have an equivalence relation $\groupoid{Q}{W}$ that is described
by \etale{} fixed-point reflecting morphisms $\map{e}{W}{Q}$ (the diagonal),
$\map{i}{Q}{Q}$ (switching the two factors) and $\map{c}{Q\times_W Q}{Q}$
(projection onto the first and third factors of $W\times_X W\times_X W$). As
the descent condition is satisfied, these morphisms descend to \etale{}
morphisms $\map{e}{Z_W}{Z_Q}$, $\map{i}{Z_Q}{Z_Q}$ and $\map{c}{Z_Q\times_{Z_W}
  Z_Q}{Z_Q}$ so that $\groupoid{Z_Q}{Z_W}$ is a groupoid in algebraic
spaces. It remains to verify that $Z_Q\to Z_W\times Z_W$ is a monomorphism, or
equivalently, that the stabilizer $Z_Q\times_{Z_W\times Z_W} Z_W\to Z_W$ is an
isomorphism.

It is enough to show that the open immersion
$\map{e\times\id{Z_W}}{W\times_{Z_W\times Z_W} Z_W}{Q\times_{Z_W\times Z_W}
  Z_W}$ is surjective. Let $\map{q}{W}{Z_W}$ denote the quotient map and let
$\map{(w_1,w_2)}{\Spec(k)}{Q=W\times_X W}$ be a
point such that $q\circ w_1=q\circ w_2$. After replacing
$k$ with a field extension, we can assume that $(w_1,w_2)$ is
in the image of an element $\map{r}{\Spec(k)}{R_W}$. Now consider the cartesian
diagram
$$\xymatrix{%
\stab(W)\ar@{(->}[r]\ar[d] & \stab(X)\times_X W\ar@{(->}[r]\ar[d] & R_W\ar[d]\\
W\ar@{(->}[r]^e & Q\ar@{(->}[r]\ar@{}[ul]|\square & W\times W.\ar@{}[ul]|\square
}$$
As $f$ is fixed-point reflecting, the leftmost arrow in the first row is an
isomorphism and it follows that $(w_1,w_2)$ factors through the diagonal
$\map{e}{W}{Q}$.

Let $Z_X$ be the quotient sheaf of the equivalence relation in the
\etale{} topology. This is an algebraic space and there is a canonical morphism
$X\to Z_X$ making the diagram~\eqref{E:existence-GC-fpr} cartesian.
As strongly geometric quotients and the descent condition are descended by
\etale{} base change by Proposition~\pref{P:quotients-and-base-change} and
Corollary~\pref{C:dc-and-base-change}, it follows that $X\to Z_X$ is a \GC{}
quotient.
\end{proof}

\end{section}


\begin{section}{Finite quotients of affine and \AF{}-schemes}
\label{S:affine-case}

In this section, we review the known results on quotients of finite
locally free groupoids of affine schemes. These are then easily extended to
groupoids of schemes such that every orbit is contained in an affine open
subscheme. The general existence results were announced in~\cite[No.~212]{fga}
by Grothendieck and proven in~\cite{gabriel_quotients} by Gabriel. An
exposition of these results with full proofs can also be found
in~\cite[Ch.~III, \S2]{gabriel_algebraic_groups}. 

Besides the existence results, a list of properties of the quotient when it
exists is given in Proposition~\pref{P:finite-GC-quotient:properties}. This
proposition is also valid for algebraic spaces.


\begin{theorem}[{\cite[No.~212, Thm.~5.1]{fga}}]
\label{T:existence-of-quotient:affine}
Let $S=\Spec(A)$, $X=\Spec(B)$ and $R=\Spec(C)$ be affine schemes and
let $\groupoid{R}{X}$ be a finite locally free $S$-groupoid. Let
$\map{p_1,p_2}{B}{C}$ be the homomorphisms corresponding to $s$ and $t$. Let
$Y=\Spec\left(B^R\right)$ where $B^R$ is the equalizer of the homomorphisms
$p_1$ and $p_2$. Let $\map{q}{X}{Y}$ be the morphism corresponding to the
inclusion $B^R\inj B$. Then $q$ is integral and a \GC{} quotient.
\end{theorem}
\begin{proof}
That $q$ is integral and a geometric quotient is proven
in~\cite[Thm.~4.1]{gabriel_quotients}. As $j_{/Y}$ is proper we have that $q$ is
a strongly geometric quotient. It satisfies the descent condition by
Theorem~\pref{T:descent-condition-holds-for-univ-open/int} and is thus a \GC{}
quotient by definition.
\end{proof}

\begin{remark}
Theorem~\pref{T:existence-of-quotient:affine} has a long history and there are
several proofs, e.g.~\cite[Prop.~5.1]{keel_mori_quotients}
and~\cite[\S3]{conrad_coarsespace}.
In the classical situation, the groupoid is induced by a finite
group or, more generally, a finite flat group scheme.
A good exposition of the theorem for finite
groups is given in~\cite[Exp.~V]{sga1}. For algebraic group schemes, the result
can be found in~\cite[Thm.~1, p.~111]{mumford_abel_vars}. For general group
schemes, a proof is given in~\cite[\S 4]{geer-moonen_abel-vars}. 
\end{remark}

\begin{lemma}\label{L:affine-nbhds-of-orbits}
Let $\groupoid{R}{X}$ be a finite locally free groupoid of schemes. Let
$Z\subseteq |X|$ be a finite set of points and let $U\subseteq X$ be an affine
open neighborhood of the orbit $R(Z)$. Then there exists an $R$-stable affine
open neighborhood $U'\subseteq U$ of the orbit $R(Z)$.
\end{lemma}
\begin{proof}
This is proved in~\cite[5b)]{gabriel_quotients}, based on the proof
of~\cite[Exp.~VIII, Cor.~7.6]{sga1}.
\end{proof}

\begin{theorem}[{\cite[No.~212, Thm.~5.3]{fga}}]
\label{T:existence-of-quotient:schemes}
Let $\groupoid{R}{X}$ be a finite locally free groupoid of
\emph{schemes}. Then a \GC{} quotient $\map{q}{X}{Y}$ with $q$ affine and $Y$ a
\emph{scheme} exists if and only if every $R$-orbit of $|X|$ is contained in an
affine open subscheme.
\end{theorem}
\begin{proof}
The necessity is obvious. To prove sufficiency, let $x\in |X|$. Then by
assumption and Lemma~\pref{L:affine-nbhds-of-orbits}, there is an affine open
$R$-stable neighborhood $U$ of the orbit $R(x)$. It is enough to prove
the theorem after replacing $(R,X)$ with $(R|_U,U)$ as a \GC{} quotient is
categorical. This is Theorem~\pref{T:existence-of-quotient:affine}.
%
%
\end{proof}




\begin{remark}\label{R:existence-of-quotient:scheme}
In Theorem~\pref{T:existence-of-quotient:alg_spaces}, we will show that if
$X\to S$
is a \emph{separated} algebraic space, then a \GC{} quotient $\map{q}{X}{Y}$
exists and is \emph{affine}. Thus, if $X\to S$ is a separated scheme, then it
follows, from Theorem~\pref{T:existence-of-quotient:schemes}, that a
geometric quotient $Y=X/R$ exists as a \emph{scheme} if and only if every
$R$-orbit of $|X|$ is contained in an affine open subset.
\end{remark}

\begin{remark}
When we replace the groupoid with a finite group, then
Theorem~\pref{T:existence-of-quotient:schemes} is a classical
result. It can be traced back to Serre~\cite[Ch.~III, \S12,
Prop.~19]{serre_groupes_alg_corps_de_classes} when $X$ is an algebraic
variety. Also see~\cite[Thm.~1, p.~111]{mumford_abel_vars} for the case when
$X$ is an algebraic scheme and~\cite[Exp.~V \S1]{sga1} or
\cite[Thm.~4.16]{geer-moonen_abel-vars} for arbitrary schemes.
\end{remark}

Recall that a scheme is \AF{} if every finite set of points has an open affine
neighborhood. Clearly, Theorem~\pref{T:existence-of-quotient:schemes} applies
to any \AF{}-scheme $X$. For the definition and properties of \AF{} morphisms
we refer to Appendix~\ref{A:AF-morphisms}.

\begin{proposition}\label{P:finite-GC-quotient:properties}
Let $\groupoid{R}{X}$ be a finite locally free groupoid of algebraic
spaces and assume that a geometric quotient $\map{q}{X}{Y}$ exists and that
$q$ is \emph{affine}\footnote{The assumption that $q$ is affine can be replaced
with the condition that $q$ is separated
using~\cite[Thm.~8.5]{rydh_noetherian-approx}.}.
Then $q$ is an integral and universally open \GC{}
quotient and $j_{/Y}$ is proper.
Consider the following properties of a morphism of algebraic spaces:
\begin{enumerate}
\myitem{A} quasi-compact, universally closed, universally open,\\
separated, quasi-separated, affine, quasi-affine, \AF{}; \label{I:FGC-A}
\myitem{B} finite type, locally of finite type, proper; \label{I:FGC-B}
\myitem{B$'$} projective, quasi-projective. \label{I:FGC-B'}
\end{enumerate}
If $X\to S$ has one of the properties in \ref{I:FGC-A}, then $Y\to S$ has the
corresponding property. The same holds for the properties in \ref{I:FGC-B} if
$S$ is
locally noetherian and for those
in \ref{I:FGC-B'} if $S$ is noetherian.
\end{proposition}
\begin{proof}
  As $s$ is universally open and proper, it follows by
  Proposition~\pref{P:quotients-misc-props} that $q$ is universally open and
  universally closed. In particular $q$ is integral~\cite[Prop.~18.2.8]{egaIV}.
  As $j_{/Y}$ is proper we have that $q$ is strongly
  topological. Furthermore, $q$ satisfies the descent condition by
  Theorem~\pref{T:descent-condition-holds-for-univ-open/int} and is thus a \GC{}
  quotient.

The statement about the first three properties in \ref{I:FGC-A} follows
immediately as $q$ is surjective.
The properties ``separated'' and ``quasi-separated'' follow from
Propositions~\pref{P:proper-actions}
and~\pref{P:proper-actions-for-proper-groups}. For properties ``affine'' and
``quasi-affine'' we can assume that $S$ is affine and then the
property ``affine'' follows from Theorem~\pref{T:existence-of-quotient:affine}.

Assume that $X$ is quasi-affine. To show that $Y$ is quasi-affine it is
enough to show that there is an affine open covering of the form $\{Y_f\}$
with $f\in\Gamma(Y)$. Let $y\in Y$ be a point and let $q^{-1}(y)$ be the
corresponding
orbit in $X$. Then, as $X$ is quasi-affine, there is a global section $g\in
\Gamma(X)$ such that $X_g$ is an affine neighborhood of $q^{-1}(y)$,
cf.~\cite[Cor.~4.5.4]{egaII}. Let $f=\norm_t(s^*g)\in \Gamma(X)$ be the norm
of $s^*g$ along $t$~\cite[6.4.8]{egaII}. Then $f$ is invariant,
i.e., $s^*f=t^*f$, and $X_f\subseteq X_g$ is an affine neighborhood of
$q^{-1}(y)$, cf.\ the proof of~\cite[Exp.~VIII, Cor.~7.6]{sga1}. Thus $f\in
\Gamma(Y)$ and $q(X_f)=Y_f$ is a geometric quotient of the groupoid
$\groupoid{R_f}{X_f}$. As $X_f$ is affine so is $Y_f$ by \ref{I:FGC-A}
for ``affine''. This shows that $Y$ is quasi-affine.

Finally, assume that $X\to S$ is \AF{}. Let $S'$ be an affine scheme and let
$S'\to S$ be a morphism so that $X'=X\times_S S'$ is an \AF{}-scheme. Let
$Y'=Y\times_S S'$ and $R'=R\times_S S'$ so that $\map{q'}{X'}{Y'}$ is a
topological quotient of $\groupoid{R'}{X'}$. We will show that $Y'$ is an
\AF{}-scheme. Let $Z\subseteq |Y'|$ be a finite subset. Then $q'^{-1}(Z)$ is a
finite subset of $X'$ and therefore admits an $R'$-stable affine open
neighborhood $U'\subseteq X'$ according to
Lemma~\pref{L:affine-nbhds-of-orbits}. By
Theorem~\pref{T:existence-of-quotient:affine}, there is a \GC{} quotient $U'\to
V'$ where $V'$ is an affine scheme. As $q'$ is a topological quotient, the
induced morphism $\map{r}{V'}{q'(U')}$ is a separated universal homeomorphism.
Then $r$ is integral~\cite[Cor.~5.22]{rydh_submersion_and_descent} and it
follows that $q'(U')$ is affine by Chevalley's
theorem~\cite[Thm.~8.1]{rydh_noetherian-approx}.

Now assume that $S$ is locally noetherian. As we have already shown the
statement for
quasi-compact, universally closed and separated, it is enough to show the
statement for the property ``locally of finite type'' in group \ref{I:FGC-B}.
Assume that $X\to S$ is locally of finite type. Then $q$ is finite. As the
quotient is uniform, we can, in order to prove that $Y\to S$ is locally of
finite type, assume that $S$ and $Y$ are affine and hence also $X$. It is then
easily
seen that $Y\to S$ is of finite type. For details, see the argument
in~\cite[Exp.~V, Cor.~1.5]{sga1}.

For the properties in \ref{I:FGC-B'}, we cannot assume that $S$ is affine as
projectivity and quasi-projectivity are not local on the base. The statement
about (quasi\nobreakdash-)projectivity is probably well-known but I am not aware
of any full proof. A sketch is given
in~\cite[Ch.~IV, Prop.~1.5]{knutson_alg_spaces} and is also discussed in
\cite{rydh_gammasymchow_inprep}, but these proofs are for quotients by
finite groups. We will now prove the general case.

Let $\sL$ be an ample
sheaf on $X$ and let
$\sL'=\norm_t(s^*\sL)=\norm_{t_*\sO_R/\sO_X}(t_*s^*\sL)$. This is an ample
invertible sheaf by~\cite[Cor.~6.6.2]{egaII}. Moreover, it comes with a
canonical $R$-linearization~\cite[Ch.~1, \S3]{GIT}, i.e., a canonical
isomorphism $\map{\phi}{s^*\sL'}{t^*\sL'}$ satisfying a cocycle
condition. This is obvious from the description $\sL'=p^*\norm_{X/\stX}(\sL)$
where $\map{p}{X}{\stX}$ is the stack quotient of
$\groupoid{R}{X}$. Consider the graded $\sO_X$-algebra $\sA=\oplus_{n\geq
0} \sL'^n$. As $\sL'$ is ample, we have a canonical (closed) immersion $X\inj
\Proj(f_*\sA)$ where $\map{f}{X}{S}$ is the structure morphism.

Let $(f_*\sA)^R$ be
the invariant ring, where $(f_*\sL'^n)^R$ is the equalizer of
$$\xymatrix{%
f_*\sL'^n\ar[r]^-{s^*}\ar@/_1pc/[rr]_{t^*} &
(f\circ s)_*s^*\sL'^n\ar[r]^-{\phi} &
(f\circ s)_*t^*\sL'^n.}$$
It can then be shown, analogous to the
case of a finite group action, that the quotient $Y$ is a subscheme of
$\Proj_S((f_*\sA)^R)$. As $S$ is locally noetherian, it follows that
$(f_*\sA)^R$ is
a finitely generated $\sO_S$-algebra, but it is not necessarily generated by
elements of degree $1$. As $S$ is noetherian there is an integer $m$ such that
$(f_*\sA^R)^{(m)}=\oplus_{n\geq 0} f_*(\sL'^{mn})^R$ is generated in degree
$1$. Hence $Y$ is (quasi\nobreakdash-)projective.
\end{proof}

\begin{remark}
If $S$ is of characteristic zero, i.e., a $\Q$-space, then the \GC{} quotient
$\map{q}{X}{Y}$ of the proposition is \emph{universal}, i.e., it commutes with
any base change. In fact, there exists a Reynolds operator, i.e., an
$\sO_Y$-module retraction $\map{\cR}{q_*\sO_X}{\sO_Y}$ of the inclusion
$\sO_Y=(q_*\sO_X)^R \inj q_*\sO_X$. The sequence~\eqref{E:geometric-quot-seq}
is thus split exact which shows that $q$ is a universal geometric quotient. We
construct the Reynolds operator as follows: the rank $r$ of $s$ is constant on
each connected component and constant on $R$-orbits. Locally the Reynolds
operator $\cR$ is defined by $\frac{1}{r}\trace_{s^*}\circ t^*$ where
$\map{s^*,t^*}{q_*\sO_X}{(q\circ s)_*\sO_R}$ are the $\sO_Y$-homomorphisms
induced by $s$ and $t$.

%
More generally, in any characteristic, the quotient is universal if the stack
$[\groupoid{R}{X}]$ is
\emph{tame}~\cite{abramovich_olsson_vistoli_tame_stacks}.
\end{remark}

\begin{remark}\label{R:existence-of-quotient:AF}
When $S$ is a \emph{scheme} and $X\to S$ is \AF{}, then any $R$-orbit of $|X|$
is
contained in an affine open subscheme and the conclusion of
Theorem~\pref{T:existence-of-quotient:schemes} holds. In particular, this is
true for $X\to S$ (quasi\nobreakdash-)affine or
(quasi\nobreakdash-)projective. Furthermore,
Proposition~\pref{P:finite-GC-quotient:properties} shows that geometric
quotients exist in the following categories:
\begin{enumerate}\renewcommand{\theenumi}{(\arabic{enumi})}
\item schemes affine over $S$;
\item schemes quasi-affine over $S$;
\item schemes that are \AF{} over $S$;
\item schemes that are projective over a noetherian base scheme $S$; and
\item schemes that are quasi-projective over a noetherian base scheme $S$.
\end{enumerate}
\end{remark}

\end{section}


\begin{section}{Finite quotients of algebraic spaces}
\label{S:finite-quotients}

Let $\groupoid{R_X}{X}$ be a groupoid. For any \etale{} morphism $U\to X$
we will construct a groupoid $\groupoid{R_W}{W}$ with a square \etale{}
morphism $\map{h}{(R_W,W)}{(R_X,X)}$. The construction requires that
$\groupoid{R_X}{X}$ is \emph{proper}, flat and of finite presentation. If
$\groupoid{R_X}{X}$ is \emph{finite} and $U\to X$ is surjective, then
$h|_{\fpr}$ will be surjective. Using Theorem~\pref{T:existence-of-GC-by-fpr}
and the results of Section~\ref{S:affine-case}, we then deduce the existence of
a quotient $X/R$.

\begin{proposition}\label{P:cover-of-finite-groupoid}
Let $\groupoid{R_X}{X}$ be a groupoid that is proper, flat and of finite
presentation and let $\map{f}{U}{X}$ be an \etale{} and separated
morphism. Then there is a groupoid $(R_W,W)$ together with a square separated
\etale{} morphism $\map{h}{(R_W,W)}{(R_X,X)}$ and an \etale{} and
separated morphism $W\to U$. In particular, if $U$ is an \AF{}-scheme (e.g., a
disjoint union of affine schemes), then so is $W$.
If, in addition, $\groupoid{R_X}{X}$ is \emph{finite} and $f$ is surjective,
then $h|_{\fpr(h)}$ is surjective.
\end{proposition}
\begin{proof}
We will construct $W$ functorially such that a point in the
fiber $W_x$ over a point $x\in |X|$ corresponds to the choice of a point in the
fiber $U_{x'}$ for every point $x'$ in the orbit of $x$. More precisely,
given an $X$-scheme $T$, an $X$-morphism $T\to W$ corresponds to a section of
$$\map{\pi_{12}}{T\times_{g,X,s} R_X \times_{t,X,f} U}{ T\times_{g,X,s} R_X}$$
where $\map{g}{T}{X}$ is the structure morphism. Thus, $W$ is the
\emph{Weil restriction} $\weilr_{s}(R_X\times_{t,X,f} U/R_X)$,
cf.\ Appendix~\ref{A:Weil-restrict}. As $\map{s}{R_X}{X}$ is proper, flat and
of finite
presentation, the functor $W$ is an algebraic space, that is
separated and locally of finite presentation over $X$.
By Proposition~\pref{P:PI-prop}, the morphism $W\to X$ is \etale{} and
separated. Furthermore, by the same proposition, the unit section of
$\map{s}{R_X}{X}$ gives rise to a
factorization of $W\to X$ into an \etale{} and separated morphism
$W\to U$ followed by $f$.
If $U$ is an \AF{}-scheme, then so is $W$ since $W\to U$ is \AF{}.

We obtain an easier description of $W$ using the stack
$\stX=[\groupoid{R_X}{X}]$. Then $W=\stW\times_\stX X$ where
$\stW=\weilr_{X/\stX}(U/X)$.
This induces a groupoid $(R_W,W)$ with $R_W=W\times_{\stW} W$ and the morphism
$\stW\to \stX$ induces a square morphism $(R_W,W)\to (R_X,X)$. The morphism $W\to
U$ is given by the adjunction formula
$\Hom_X(T\times_{\stX} X, U)=\Hom_{\stX}(T,\stW)$ with $T=\stW$.


\newcommand{\kbar}{\overline{k}}

Finally, we show that $h|_{\fpr(h)}$ is surjective when $R_X\to X$ is finite
and $\map{f}{U}{X}$ is surjective. Let
$\map{x}{\Spec(\kbar)}{X}$ be a geometric point of $X$. A lifting
$\map{w}{\Spec(\kbar)}{W}$ of $x$ corresponds to a morphism
$\map{\varphi}{s^{-1}(x)}{U}$ such that $t=f\circ\varphi$. Let
$R_x=s^{-1}(x)_\red$
which we consider as an $X$-scheme using $t$. As $U\to X$ is \etale{}, any
$X$-morphism $R_x\to U$ induces a unique morphism $\varphi$ as above. If
$R_X\to X$ is finite, then $R_x$ is a finite set of points. We may then choose
an $X$-morphism $R_x\to U$ such that its image contains at most one point in
every
fiber of $\map{f}{U}{X}$. This corresponds to a point $w$ in the fixed-point
reflecting locus of $W\to X$.
\end{proof}

\begin{remark}
Let $G=\{g_1,g_2,\dots,g_n\}$ be a finite group acting on an
algebraic space $X$ and let $\groupoid{R_X=G\times X}{X}$ be the induced
groupoid. Let
$\map{f}{U}{X}$ be an \etale{} and separated morphism. Then the \etale{} cover
$W\to X$ of Proposition~\pref{P:cover-of-finite-groupoid} is the fiber product
of $g_1\circ f,g_2\circ f,\dots,g_n\circ f$. The morphism $W\to U$ is the
projection on the factor corresponding to the identity element $g_i=e\in G$.
\end{remark}

\begin{theorem}\label{T:existence-of-quotient:alg_spaces}
Let $\groupoid{R}{X}$ be a finite locally free groupoid with finite
stabilizer $\stab(X)=j^{-1}\bigl(\Delta(X)\bigr)\to X$. Then a \GC{} quotient
$\map{q}{X}{X/R}$ exists and $q$ is affine. Hence it has the properties of
Proposition~\pref{P:finite-GC-quotient:properties}.
\end{theorem}
\begin{proof}
The question is \etale{}-local over $S$ so we can assume that $S$ is
affine. Let $\map{\varphi}{U}{X}$ be an \etale{} presentation such that $U$ is a
disjoint union of affine schemes. Let $\map{h}{W}{X}$ be the \etale{} cover
constructed in Proposition~\pref{P:cover-of-finite-groupoid} using $U$. As the
stabilizer
is finite, the subset $\fpr(h)\subseteq |W|$ is open by
Proposition~\pref{P:fpr-of-unramified-is-open}. Furthermore, $W|_\fpr\to X$ is
surjective by Proposition~\pref{P:cover-of-finite-groupoid}. Thus $W|_\fpr\to
X$ is an \etale{} square fpr cover such that $W|_\fpr$ is an \AF{}-scheme. By
Theorem~\pref{T:existence-of-quotient:schemes} and
Remark~\pref{R:existence-of-quotient:AF}, a \GC{} quotient of $W|_\fpr$ exists.
Hence a geometric quotient $X/R$ exists by
Theorem~\pref{T:existence-of-GC-by-fpr}.
\end{proof}

\begin{corollary}[Deligne]\label{C:existence-of-group-quotient:alg_spaces}
Let $G\to S$ be a finite locally free group scheme acting on a separated
algebraic space $X\to S$. Then a \GC{} quotient $\map{q}{X}{X/G}$ exists, $q$
is affine and $X/G\to S$ is separated.
\end{corollary}
\begin{proof}
As $X$ is separated, the finite locally free groupoid
$\groupoid{G\times_S X}{X}$ has finite diagonal. In
particular, its stabilizer is finite. The existence of a \GC{} quotient $q$
thus follows from Theorem~\pref{T:existence-of-quotient:alg_spaces}.
\end{proof}


For symmetric products we can find more explicit \etale{} covers as follows.

\begin{corollary}\label{C:existence-of-Sym:sep-alg-space}
Let $X\to S$ be a separated morphism of algebraic spaces. Then a
\GC{} quotient $\Sym^d(X/S):=(X/S)^d/\SG{d}$ exists as a separated algebraic
space. Let $\{S_\alpha\to S\}_\alpha$ and
$\{U_\alpha\to X\times_S S_\alpha\}_\alpha$ be sets of \etale{} morphisms of
separated algebraic spaces. Then the diagram
\begin{equation}\label{E:Sym-diagram}
\vcenter{\xymatrix{
\coprod_\alpha (U_\alpha/S_\alpha)^d|_\fpr\ar[r]\ar[d] & (X/S)^d\ar[d] \\
\coprod_\alpha \Sym^d(U_\alpha/S_\alpha)|_\fpr\ar[r]
  & \Sym^d(X/S)\ar@{}[ul]|\square
}}
\end{equation}
is cartesian and the horizontal morphisms are \etale{}. If the $U_\alpha$'s are
such that for every $s\in |S|$, any set of $d$ points in $X_s$ lies in
the image of some $U_\alpha$, then the horizontal morphisms are surjective.
In particular, there is an \etale{} cover of $\Sym^d(X/S)$ of the form
$\left\{\Sym^d(U_\alpha/S_\alpha)|_{V_\alpha}\right\}$ with $U_\alpha$ and
$S_\alpha$ affine and $V_\alpha$ an open subset.
\end{corollary}
\begin{proof}
Corollary~\pref{C:existence-of-group-quotient:alg_spaces} shows the
existence of the \GC{} quotient $\Sym^d(X/S)$. As $(U_\alpha/S_\alpha)^d\to
(X/S)^d$ is \etale{}, the diagram~\eqref{E:Sym-diagram} is cartesian by the
descent condition. Let $\map{x=(x_1,x_2,\dots,x_d)}{\Spec(k)}{(X/S)^d}$
be a geometric point. If
$x_1,x_2,\dots,x_d$ lies in the image of some $U_\alpha$ we can choose a
lifting $\map{u=(u_1,u_2,\dots,u_d)}{\Spec(k)}{(U_\alpha/S_\alpha)^d}$
such that $u_i=u_j$ if and only if $x_i=x_j$. Then $u$ is a point in
the fixed-point reflecting locus of $(U_\alpha/S_\alpha)^d\to (X/S)^d$. This
shows the surjectivity of the horizontal morphisms.
\end{proof}

\begin{remark}
Note that the stabilizer of the permutation action of $\SG{d}$ on $(X/S)^d$
is proper exactly when $X\to S$ is separated. Thus,
the proof of Corollary~\pref{C:existence-of-Sym:sep-alg-space} definitively
fails unless $X\to S$ is separated.
\end{remark}


\end{section}


\begin{section}{Coarse moduli spaces of stacks}
\label{S:cms-stacks}
In this section, we will prove the existence of a coarse moduli space to any
algebraic stack $\stX$ with finite inertia. We proceed as follows.
\begin{enumerate}
\item\label{E:step-1}
We find a quasi-finite flat presentation of
$\stX$~\cite[Thm.~7.1]{rydh_etale-devissage}.
\item\label{E:step-2}
We find an \etale{} representable cover $\stW\to \stX$ such that there
exists an \AF{}-scheme $V$ and a finite flat presentation $V\to \stW$
(Proposition~\ref{P:cover-of-qf-stack}).
\item\label{E:step-3}
We show that $\stW\to \stX$ is fixed-point reflecting over an open substack
$\stW|_{\fpr}\subseteq \stW$ and that $\stW|_{\fpr}\to \stX$ is surjective.
(Proposition~\ref{P:fpr-of-unramified-is-open:stacks}).
\item\label{E:step-4}
We deduce the existence of a coarse moduli space to $\stX$ from the
existence of a coarse moduli space to $\stW$
(Theorem~\ref{T:existence-of-GC-by-fpr:stacks}).
\end{enumerate}
The assumption that $\stX$ has finite inertia, and not merely quasi-finite
inertia, is only used in
step~\ref{E:step-3}.

Keel and Mori~\cite{keel_mori_quotients} more or less proceed in the same
way. Using stacks, as in~\cite{conrad_coarsespace}, instead of
groupoids, as in~\cite{keel_mori_quotients}, gives a more streamlined
presentation and simplifies many arguments. In particular, the reduction
from the quasi-finite case to the finite case becomes much more
transparent. In step~\ref{E:step-4}, we use the descent condition
(Definition~\ref{D:descent-condition:stacks}), which simplifies earlier proofs,
cf.~\cite[Thm.~3.1 and
Thm.~4.2]{conrad_coarsespace}. We also avoid the somewhat complicated limit
methods used in~\cite[\S5]{conrad_coarsespace} and obtain a slightly more
general result.

We begin by rephrasing the results of \S\ref{S:fpr-and-descent} in stack
language.
If $U\to \stX$ is a flat presentation and $R=U\times_{\stX} U$, then there is a
one-to-one correspondence between equivariant morphisms $U\to Z$ with respect to
the groupoid
$\groupoid{R}{U}$ and morphisms $\stX\to Z$. We say that a morphism
$\stX\to Z$ is a topological (resp.\ geometric, resp.\ categorical,\ \dots)\
quotient if $U\to \stX\to Z$ is such a quotient. This definition does not depend
on the choice of presentation $U\to \stX$ and can be rephrased as follows.

\begin{definition}\label{D:quotient-of-stacks}
Let $\stX$ be an algebraic stack, let $Z$ be an algebraic space and let
$\map{q}{\stX}{Z}$ be
a morphism. Then $q$ is:
\begin{enumerate}
\item \emph{categorical} if $q$ is initial among morphisms from $\stX$ to
algebraic spaces;
\item \emph{topological} if $q$ is a universal homeomorphism;
\item \emph{strongly topological} if $q$ is a strong homeomorphism,
cf.\ Appendix~\ref{A:strong-homeomorphisms};
\item \emph{geometric} if $q$ is a universal homeomorphism and
$\sO_Z\to q_*\sO_{\stX}$
is an isomorphism; and
\item \emph{strongly geometric} if $q$ is a strong homeomorphism and
$\sO_Z\to q_*\sO_{\stX}$ is an isomorphism.
\end{enumerate}
\end{definition}

\begin{remark}
If $\map{q'}{\stX}{Z'}$ and $\map{r}{Z'}{Z}$ are universal homeomorphisms (so
that $q'$ and $q=r\circ q'$ are topological quotients), then $r$ need not be
separated. If $q$ is a strong homeomorphism (i.e., $q$ is a strongly
topological quotient), then $r$ is necessarily separated by
Corollary~\pref{C:strong-homeo-and-separatedness}. Thus, a strongly topological
quotient $\map{q}{\stX}{Z}$ is ``maximally non-separated'' among topological
quotients.
\end{remark}

Let $\map{f}{\stX}{\stY}$ be a morphism of stacks. There is then an induced
morphism $\map{\varphi}{I_{\stX}}{I_{\stY}\times_{\stY} \stX}$. If
$\map{x}{\Spec(k)}{\stX}$ is a point and $y=f\circ x$, then $\varphi_x$ is the
natural morphism of $k$-groups $\Isom_\stX(x,x)\to \Isom_\stY(y,y)$.


\begin{definition}
Let $\map{f}{\stX}{\stY}$ be a morphism of stacks. We say that $f$ is
\emph{fixed-point reflecting} (or stabilizer preserving), abbreviated fpr, at
a point $x\in |\stX|$ if
$\map{\varphi_x}{I_{\stX}\times_{\stX} \Spec(k)}{I_{\stY}\times_{\stY}
\stX\times_{\stX} \Spec(k)}$ is an isomorphism for some representative
$\Spec(k)\to \stX$ of $x$.
We say that $f$ is fixed-point reflecting if $f$
is fixed-point reflecting at every point. We let $\fpr(f)\subseteq |\stX|$ be
the subset over which $f$ is fixed-point reflecting.
\end{definition}

It is well-known that $f$ is representable if and only if $\varphi$ is a
monomorphism. If $f$ is unramified, then $f$ is fixed-point reflecting at every
point if and only if $\varphi$ is an isomorphism.

\begin{remark}
Let $\stX$ and $\stY$ be stacks with presentations $U\to \stX$ and $V\to \stY$
and let $R_U=U\times_{\stX} U$ and $R_V=V\times_{\stY} V$.
Then there is a one-to-one correspondence between square morphisms $(R_U,U)\to
(R_V,V)$ and morphisms $\stX\to \stY$ together with an isomorphism
$U\to V\times_{\stY} \stX$. Under this correspondence, $U\to V$ is fixed-point
reflecting if and only if $\stX\to \stY$ is so.
\end{remark}

The following is a reformulation of
Proposition~\pref{P:fpr-of-unramified-is-open} for stacks.

\begin{proposition}\label{P:fpr-of-unramified-is-open:stacks}
Let $\map{f}{\stX}{\stY}$ be a representable and unramified morphism of
stacks. If the inertia stack $I_{\stY}\to \stY$ is universally closed, then the
subset
$\fpr(f)\subseteq |\stX|$ is open.
\end{proposition}

\begin{definition}\label{D:descent-condition:stacks}
Let $\stX$ be an algebraic stack and let $\map{q}{\stX}{X}$ be a topological
quotient.
We say that $q$ satisfies the \emph{descent condition} if for any \etale{}
fixed-point reflecting morphism $\map{f}{\stW}{\stX}$ there exists a
topological quotient
$\stW\to W$ and a $2$-cartesian square
\begin{equation*} 
\xymatrix{
\stW\ar[r]^f\ar[d] & \stX\ar[d]^q \\
W\ar[r] & X\ar@{}[ul]|\square
}\end{equation*}
where $W\to X$ is \emph{\etale}.
\end{definition}

\begin{proposition}\label{P:geom+dc=>cat:stacks}
A strongly geometric quotient $\stX\to X$ satisfies the descent condition
for $\Et_{\qsep}(\stX)$ and is categorical.
\end{proposition}
\begin{proof}
See Theorem~\pref{T:descent-condition-holds-for-univ-open/int} and
Proposition~\pref{P:geom+dc=>cat}.
\end{proof}

\begin{definition}
We say that $\stX\to X$ is a \emph{coarse moduli space} if $\stX\to X$ is
a strongly geometric quotient. As a strongly geometric quotient is categorical,
we will speak about \emph{the} coarse moduli space when it exists.
\end{definition}

\begin{remark}\label{R:dc-GC-and-stacks}
Let $\stX$ be an algebraic stack with a flat presentation $U\to \stX$ and let
$R=U\times_{\stX} U$. Let
$\map{q}{\stX}{X}$ be a topological quotient. Then $q$ satisfies the descent
condition (resp.\ is a coarse moduli space) if and only if $U\to\stX\to X$
satisfies the descent condition of Definition~\pref{D:descent-condition}
(resp.\ is a \GC{} quotient) for the groupoid $(R,U)$.
\end{remark}

\begin{theorem}\label{T:existence-of-GC-by-fpr:stacks}
Let $\map{f}{\stW}{\stX}$ be a surjective, \etale{}, separated and fpr morphism
of algebraic
stacks. Let $\stQ=\stW\times_{\stX} \stW$. If $\stW$ has a coarse moduli space
$W$,
then coarse moduli spaces $\stQ\to Q$ and $\stX\to X$ exist. Furthermore, the
diagram
\begin{equation*}
\xymatrix{\stQ \ar@<.5ex>[r] \ar@<-.5ex>[r]\ar[d] & \stW \ar[r]^f\ar[d] & \stX\ar[d] \\
Q \ar@<.5ex>[r] \ar@<-.5ex>[r] & W \ar[r] & X}
\end{equation*}
is cartesian.
\end{theorem}
\begin{proof}
This follows immediately from Theorem~\pref{T:existence-of-GC-by-fpr} and
Remark~\pref{R:dc-GC-and-stacks}.
\end{proof}

\begin{proposition}[{\cite[\S4]{keel_mori_quotients},
\cite[Lem.~2.2]{conrad_coarsespace},
\cite[Thm.~6.3]{rydh_etale-devissage}}]\label{P:cover-of-qf-stack}
Let $\stX$ be an algebraic stack with locally quasi-finite and separated
diagonal.
Then there is a representable
\etale{} separated morphism $\map{h}{\stW}{\stX}$ such that $\stW$ has a
\emph{finite}
flat presentation $V\to \stW$ with $V$ an \AF{}-scheme.
Moreover, $\fpr(h)\subseteq |\stW|\to |\stX|$ is surjective.
\end{proposition}
\begin{proof}
By~\cite[Thm.~7.1]{rydh_etale-devissage} (or~\cite[Lem.~04N0]{stacks-project} if
$\Delta_{\stX}$ is not quasi-compact) there is a locally quasi-finite flat
presentation $\map{p}{U}{\stX}$
with $U$ a scheme. Taking an open covering, we can assume that $U$ is a
disjoint union of affine schemes, hence \AF{}.
Let $\stH=\Hilb(U/\stX)\to \stX$ be the
Hilbert functor. As $U\to \stX$ is locally of finite presentation and
separated, it follows by fppf descent and the usual algebraicity result
for Hilbert functors of algebraic
spaces~\cite[Cor.~6.2]{artin_alg_formal_moduli_I}
that $\stH\to\stX$ is representable, locally of finite
presentation and separated.

As $U\to\stX$ is locally quasi-finite, a morphism $T\to \stH$
corresponds to a closed subscheme $\injmap{g}{Z}{U\times_{\stX} T}$ such that
the morphism $\map{f=\pi_2\circ g}{Z}{T}$ is
finite, flat and of finite presentation over $T$. Note that $g$ is an open
immersion if and only if $g$ is
\etale{}. As $f$ is flat and $f$ and $\map{\pi_2}{U\times_{\stX} T}{T}$ are
both locally of finite presentation, we have that the morphism $g$
is \etale{} at $z$ if and only if the fiber $g_{f(z)}$ is \etale{} at $z$
by~\cite[Rem.~17.8.3]{egaIV}. Let $Z_{\metale}$ be the open subset of $Z$ where
$g$ is \etale{}. Then the open subset $T\setminus f(Z\setminus
Z_{\metale})\subseteq T$ is the set of $t\in |T|$ such that $g_t$ is open and
closed. It follows that the substack $\stW\subset \stH$, parameterizing open and
closed subschemes, is an open substack. It is also immediately clear that
$\stW\to \stX$ is formally \etale{} and hence \etale{}.
In fact, $\stW$ can also be described as the \etale{} sheaf of pointed sets
$p_!\underline{\{0,1\}}_U$.

Let $V$ be the universal family over $\stW$. Then $V$ is an \AF{}-scheme.
In fact, as $V\inj U\times_{\stX} \stW$ is an open and closed immersion
and $\stW\to \stX$ is \etale{} and separated, we have that $V\to U$ is \AF{} so
that $V$ is an \AF{}-scheme, cf.\ Appendix~\ref{A:AF-morphisms}.

Moreover, we have a decomposition $\stW=\coprod_{d\geq 0} \stW_d$ into open and
closed substacks where $\stW_d$ parameterizes open and closed subschemes of
rank $d$ over the base. After replacing $\stW$ with $\stW\setminus \stW_0$ we
have that $V\to \stW$ is surjective.

To see that $\stW|_\fpr\to \stX$ is surjective, let $\map{x}{\Spec(k)}{\stX}$
be a point. The stabilizer group scheme $\Isom_{\stX}(x,x)$ acts on the fiber
$U_x$. If $Z\subseteq U_x$ is an open subscheme that is stable under
this action and finite over $\Spec(k)$, then the corresponding point
$\map{[Z]}{\Spec(k)}{\stW}$ is in the fixed-point reflecting locus of $\stW\to
\stX$. As $U_x$ is discrete and $\Isom_{\stX}(x,x)$ is finite, we can choose
$Z$ as any orbit of $|U_x|$.
\end{proof}

We are now ready to prove the full generalization of Keel and Mori's
theorem~\cite{keel_mori_quotients}.

\begin{theorem}\label{T:KM-general}
Let $\stX$ be an algebraic stack with finite inertia stack. Then $\stX$ has a
coarse moduli
space $\map{q}{\stX}{X}$ and $q$ is a \emph{separated} universal homeomorphism.
Let $S$ be an algebraic space and let $\stX\to S$ be a morphism.
If $\stX\to S$ is locally of finite type, then $q$ is
proper and quasi-finite.
Consider the following properties:
\begin{enumerate}
\myitem{A}\label{TI:KM:A} quasi-compact, universally closed, universally open,
separated, quasi-separated;
\label{I:SGC-A}
\myitem{B}\label{TI:KM:B} finite type, locally of finite type, proper.
\label{I:SGC-B}
\end{enumerate}
If $\stX\to S$ has one of the properties in \ref{I:SGC-A}, then $X\to S$ has the
corresponding property. If $S$ is locally noetherian, the same holds for the
properties
in \ref{I:SGC-B}.
\end{theorem}
\begin{proof}
By Propositions~\pref{P:fpr-of-unramified-is-open:stacks}
and~\pref{P:cover-of-qf-stack}, there is a representable, \etale{},
fixed-point reflecting and surjective morphism $\stW\to \stX$ and a finite flat
presentation $V\to \stW$ with $V$ an \AF{}-scheme. A
coarse moduli space $\stW\to W$ exists by
Theorem~\pref{T:existence-of-quotient:schemes}. It then follows, from
Theorem~\pref{T:existence-of-GC-by-fpr:stacks}, that a coarse moduli
space
$\map{q}{\stX}{X}$ exists, that the morphism $W\to X$ is \etale{} and that
$\stW=\stX\times_X W$.

As $\stW\to W$ is separated, so is $q$. If $\stX\to S$ is locally of finite
type, then $\stX\to X$ is locally of finite type and hence proper.
We also then have that $V\to X$ is locally quasi-finite so that $\stX\to X$ is
quasi-finite.

Among the properties in \ref{TI:KM:A}, ``separated'' and ``quasi-separated'' follow
from
Proposition~\pref{P:proper-actions} and the rest are
obvious. In \ref{TI:KM:B}, we only need to prove that if $S$ is locally
noetherian and
$\stX\to S$ is locally of finite type, then so is $X\to S$. As $\stW\to S$ is
locally of finite type, then so is $W\to S$ by
Proposition~\pref{P:finite-GC-quotient:properties}. As $W\to X$ is \etale{} and
surjective, it follows that $X\to S$ is locally of finite type.
\end{proof}


\begin{remark}
If $\stX$ is an algebraic stack such that the inertia stack $I_{\stX}\to \stX$
is flat and locally of finite presentation, then the fppf sheaf $X$ associated
to the
stack $\stX$ is a coarse moduli space of $\stX$. Indeed,
$X$ is an algebraic space and $\stX\to X$ is an fppf
gerbe~\cite[Cor.~10.8]{laumon}. It follows that $\stX\to X$ is strongly
geometric and that the formation of $X$ commutes with arbitrary base
change.

Thus, if $\stX$ has finite or flat inertia, then it has a coarse moduli space.
It is then easily seen that for the existence of a coarse moduli space, it is
not necessary that $\stX$ has proper inertia. In fact, if $X$ is any algebraic
space, $U\subseteq X$ is an open subset and $G$ is a finite group, then there
is an
algebraic stack $\stX$ with coarse moduli space $X$, stabilizer group $G$
over $U$ and trivial stabilizer group over $X\setminus U$. It is clear that
the inertia stack $I_{\stX}\to \stX$ is flat but not proper unless $U$ is also
closed.
\end{remark}

The following example shows that if $\stX$ is an algebraic stack with
quasi-finite inertia, which is neither proper nor flat,
then a coarse moduli space need not
exist. In fact, the example does not admit a topological quotient, nor a
categorical quotient. The proof hinges on the following lemma which is closely
related to Artin's example of a non-algebraic
functor~\cite[Ex.~5.11]{artin_implicit_function_thm}.

\begin{lemma}\label{L:monomorphism-affine-plane}
Let $\A{2}=\Spec(k[x,y])$ be the affine plane, let $U=D(x)$ be the complement
of the $y$-axis and let $Z=V(y)$ be the $x$-axis. Let $X$ be an algebraic space
and let $\injmap{i}{X}{\A{2}}$ be a monomorphism, locally of finite type, such
that $i|_U$ and $i|_Z$ are isomorphisms. Then there is an open neighborhood
$V\subseteq \A{2}$ of $|U|\cup |Z|$ such that $i|_V$ is an isomorphism.
\end{lemma}
\begin{proof}
We begin by noting that $X$ is a scheme~\cite[Thm.~A.2]{laumon} and that we
can replace $X$ with a quasi-compact open neighborhood of the quasi-compact
subset $|U|\cup |Z| \subseteq |X|$ so that $i$ becomes quasi-compact. Then
$i$ is quasi-finite and the lemma readily follows from Zariski's Main Theorem.
\end{proof}

\begin{example}
Let $k$ be a field of characteristic different from $2$ and let
$S=\Spec(k[x,y^2])$ be the affine plane. Let
$U=\Spec(k[x,y])$ be another affine plane, seen as a ramified double
covering of $S$. Let $\map{\tau}{U}{U}$ be the $S$-involution on $U$ given by
$y\mapsto -y$. Finally, let $W$ be the union of two copies of $U$ glued along
$x\neq 0$. On $W$ we have an involution given by simultaneously applying $\tau$
and interchanging the two copies of $U$. This gives an action of $\Z/2\Z$ on
$W$ and we let $\stX=[W/(\Z/2\Z)]$.

For us, it will more convenient to study the presentation $U\to \stX$ given by
one of the copies of $U$ in $W$. Then $U\times_{\stX} U=G\times_S U$ where $G$
is the group scheme $S\amalg (S\setminus\{x=0\})\subseteq (\Z/2\Z)_S$ and
the action of a non-trivial section of $G$ on $U$ is given by $\tau$.

Let $S'=S\setminus\{y=0\}$ and let $X'=\stX\times_S S'$.
Note that the stack $X'$ has trivial inertia, so it is an algebraic space.
The morphism $X'\to S'$ is \etale{} but not separated.
It is an isomorphism outside $\{x=0\}$ and restricted to $x=0$ it coincides
with the \etale{} double cover $U|_{S'}\to S'$.

Let $Q$ be an algebraic space over $S$ and let $\stX\to Q$ be an $S$-morphism.
Then there is an induced factorization
$${G\times_S U}\to {U\times_Q U}\inj {U\times_S U}.$$
We have that
$${U\times_S U}=\Spec\bigl(k[x,y_1,y_2]/(y_1^2-y_2^2)\bigr),$$
so that ${U\times_S U}$ is the union of two affine planes
$U_i=\Spec(k[x,t_i])$ glued along the lines $t_i=0$. In coordinates, we have
that $t_1=y_1+y_2$ and $t_2=y_1-y_2$. The image of ${G\times_S U}\to
{U\times_S U}$ is the union of $U_1$ and $U_2\setminus{\{x=0\}}$.

We now restrict everything to $U_2$. Then $(G\times_S U)|_{U_2}$ is the
disjoint union of the open subscheme $D(x)$ and the closed subscheme
$V(t_2)$. This observation, combined with
Lemma~\pref{L:monomorphism-affine-plane}, shows that
$G\times_S U\to U\times_Q U$
is not surjective. Thus, the stack $\stX$ has \emph{no topological quotient}.

In addition, the stack $\stX$ has no \emph{categorical quotient}. In fact, for
any closed point $s\in S$ on the $y^2$-axis but not on the $x$-axis, let
$Q_s\to S$ be the non-separated algebraic space which is isomorphic to $S$
outside $s$ but an \etale{} extension of degree $2$ at $s$. To be precise, over
$S'$ the space $Q_s$ is the quotient of $U|_{S'}$ by the group $S'\amalg
(S'\setminus S'_s)$ where the second component acts by $\tau$. If $k$ is
algebraically closed, then $Q_s$ is even a scheme---the affine plane with a
double point at $s$. It is clear that $\stX\to S$ factors canonically through
$Q_s$.

If a categorical quotient $\stX\to Q$ existed, then, by definition, we would
have
morphisms $Q\to Q_s$ for every $s$ as above. This shows that $U\times_Q U\inj
U\times_{Q_s} U\inj U\times_S U$ would be contained in the union of $U_1$ and
$U_2\setminus B$ where $B$ is the set of all closed points on the $t_2$-axis
except the origin. This is again impossible according to the lemma.
\end{example}

\end{section}

\appendix

\begin{section}{Descent of \'etale morphisms}\label{A:descent}
In this appendix, we review and combine the descent results for \'etale
morphisms found in~\cite[Exp.~IX]{sga1}, \cite[Exp.~VIII, \S9]{sga4}
and~\cite[\S5]{rydh_submersion_and_descent}. For an algebraic space $S$, we let
$\Et(S)=\Et_{\all}(S)$ denote the category of all \'etale morphisms $X\to S$ of
algebraic spaces. We also consider the following full subcategories of
$\Et(S)$:
\begin{align*}
\Et_{\qc}(S) &= \{ \text{\etale{} and quasi-compact morphisms} \} \\
\Et_{\qsep}(S) &= \{ \text{\etale{} and quasi-separated morphisms} \} \\
\Et_{\cons}(S) &= \{ \text{\etale{}, quasi-compact and quasi-separated
  morphisms} \} \\
\Et_{\sep}(S) &= \{ \text{\etale{} and separated morphisms} \} \\
\Et_{\sep+\qc}(S) &= \{ \text{\etale{}, separated and
quasi-compact morphisms} \} \\
\Et_{\finite}(S) &= \{ \text{\etale{} and finite morphisms} \}.
\end{align*}
Note that $\Et(S)$ is equivalent to the category of sheaves on the small
\'etale site of $S$ and that under this equivalence, $\Et_{\cons}(S)$ is
identified with the category of constructible sheaves, cf.\ \cite[Ch.~V,
  Thm.~1.5]{milne_etale_coh} or~\cite[Ch.~VII, \S1]{artin_theoremes_de_repr}.

For any morphism of algebraic spaces $\map{f}{S'}{S}$ there is a pull-back
functor $\map{f^*}{\Et(S)}{\Et(S')}$ and this makes
$\map{\Et(-)}{\AlgSp^\op}{\Cat}$ into a
pseudo-functor. By the usual theory of
descent~\cite{giraud_descente,vistoli_descent} there is also a functor
$\map{f^*_D}{\Et(S)}{\Et(S'\to S)}$ where $\Et(S'\to S)$ denotes the category of
pairs
$(X'\to S',\varphi)$ where $(X'\to S')\in\Et(S')$ and $\map{\varphi}{X'\times_S
  S'}{S'\times_S X'}$ is a descent datum.
If $(X'\to S',\varphi)$ is in the essential image of $f^*_D$, then we say that
$\varphi$ is an \emph{effective} descent datum.

For the definition of universally submersive morphisms and $f^\cons$,
see~\pref{X:cons-top} or~\cite[\S1]{rydh_submersion_and_descent}.  For the
definition of \emph{universally subtrusive} morphisms,
see~\cite[\S2]{rydh_submersion_and_descent}. Examples of universally subtrusive
morphisms are (i) morphisms that are covering in the fpqc topology,
(ii) universally open and surjective morphisms and
(iii) universally closed and surjective morphisms.
If $S$ is locally noetherian, then $\map{f}{S'}{S}$ is universally subtrusive
if and only if $f$ and $f^\cons$ are universally submersive (e.g., $f$
quasi-compact and universally submersive).

\begin{proposition}[Descent]\label{P:descent}
Let $\map{f}{S'}{S}$ be a universally submersive morphism. Then the functor
$\map{f^*_D}{\Et(S)}{\Et(S'\to S)}$ is fully faithful. Moreover, if $f^\cons$
is also
universally submersive, $(X\to S)\in \Et(S)$ and
$P\in\{\qc,\qsep,\cons,\sep,\sep+\qc,\finite\}$, then
$(X\to S)\in \Et_P(S)$ if and only if $f^*(X\to S)\in \Et_P(S')$.
\end{proposition}
\begin{proof}
The first assertion is~\cite[Exp.~VIII, Prop.~9.1]{sga4}.
The second assertion follows
from~\cite[Prop.~1.7]{rydh_submersion_and_descent},
cf.\ \cite[Prop.~5.4]{rydh_submersion_and_descent}.
\end{proof}

\begin{theorem}[Effective descent]\label{T:effective-descent}
Let $\map{f}{S'}{S}$ be a surjective morphism of algebraic spaces. Then
$\map{f^*_D}{\Et_P(S)}{\Et_P(S'\to S)}$ is an equivalence of categories in the
following
cases:
\begin{enumerate}
\myitem{1a}\label{TI:desc:prop/int} $P=\all$ and $f$ is proper or integral;
\myitem{1b}\label{TI:desc:usub+fp} $P=\all$ and $f$ is universally subtrusive
and of finite presentation;
\myitem{1c}\label{TI:desc:uopen+lfp} $P=\all$ and $f$ is universally open and
locally of finite presentation;
\myitem{1d}\label{TI:desc:fppf} $P=\all$ and $f$ is covering in the fppf
topology;
\myitem{2}\label{TI:desc:uopen} $P=\qsep$ and $f$ is universally open;
\myitem{3a}\label{TI:desc:usub+qc} $P=\cons$ and $f$ is universally subtrusive
and quasi-compact; or
\myitem{3b}\label{TI:desc:fpqc} $P=\cons$ and $f$ is covering in the fpqc
topology.
\end{enumerate}
\end{theorem}
\begin{proof}
\ref{TI:desc:prop/int} and~\ref{TI:desc:fppf} are~\cite[Exp.~VIII,
  Thm.~9.4]{sga4} and \ref{TI:desc:fppf} also follows from the more general
descent result~\cite[Cor.~10.4.2]{laumon}. Using~\ref{TI:desc:fppf} we can work
Zariski-locally and then~\ref{TI:desc:uopen+lfp} follows
from~\ref{TI:desc:usub+fp}. To prove \ref{TI:desc:usub+fp} we can assume that
$S$ is affine and then $f$ can be refined into an open covering followed by a
proper and finitely presented surjective
morphism~\cite[Thm.~3.12]{rydh_submersion_and_descent} and we conclude
by~\ref{TI:desc:prop/int} and \ref{TI:desc:fppf}. Finally~\ref{TI:desc:uopen}
and~\ref{TI:desc:usub+qc} are~\cite[Thm.~5.19 and
  Thm.~5.17]{rydh_submersion_and_descent} and~\ref{TI:desc:fpqc} is a special
case of~\ref{TI:desc:usub+qc}.
\end{proof}

I do not know if effective descent holds for all \etale{} morphisms when $f$ is
universally open or flat and quasi-compact. It is also possible that
``universally subtrusive'' can be replaced with ``universally submersive'',
cf.~\cite[Ex.~5.24]{rydh_submersion_and_descent}.

\end{section}


\begin{section}{The AF condition}\label{A:AF-morphisms}
\begin{definition}
We say that a scheme $X$ is \AF{} (affine finie) if every finite set of points
is contained in an affine open subscheme. We say that a morphism
$\map{f}{X}{S}$ of algebraic spaces is \AF{} if $X\times_S T$ is an
\AF{}-scheme for every affine scheme $T$ and morphism $T\to S$.
\end{definition}

The \AF{} condition is a natural condition for a wide range of problems. It
guarantees the existence of many objects in the category of schemes such as
finite quotients, cf.\ \cite[Exp.~V]{sga1} and
Theorem~\pref{T:existence-of-quotient:schemes},
push-outs~\cite{ferrand_conducteur} and the Hilbert scheme of
points~\cite{rydh_hilbert}. Moreover, on \AF{}-schemes \etale{}
cohomology can be calculated using \v{C}ech
cohomology~\cite[Cor.~4.2]{artin_Hensel-join},\ \cite{schroer_bigger-Brauer}.
The \AF{} condition (for $2$ points) also guarantees the existence of
embeddings into toric
varieties~\cite{wlodarczyk_toric-embeddings}.

\begin{proposition} 
Let $\map{f}{X}{Y}$ and $\map{g}{Y}{Z}$ be morphisms of algebraic
spaces.
\begin{enumerate}
\item\label{PI:AF:stability} Let $Y'\to Y$ be arbitrary. If $f$ is \AF{},
  then so is $\map{f'}{X\times_Y Y'}{Y'}$.
\item\label{PI:AF:sep} An \AF{}-scheme is separated.
\item\label{PI:AF:srepr+sep} An \AF{}-morphism is strongly representable and
  separated.
\item\label{PI:AF:total-AF} If $Y$ is \AF{}, then $f$ is \AF{} if and only if
  $X$ is \AF{}.
\item\label{PI:AF:ample} If there exists an $f$-ample invertible $\sO_X$-module,
  then $f$ is \AF{}.
  In particular, (quasi-)affine and (quasi-)projective morphisms are \AF{}.
\item\label{PI:AF:lqf+sep} If $f$ is locally quasi-finite and separated, then
  $f$ is \AF{}.
\item\label{PI:AF:comp} If $f$ and $g$ are \AF{}, then so is $g\circ f$.
\item\label{PI:AF:first} If $g\circ f$ is \AF{}, then so is $f$.
\end{enumerate}
\end{proposition}
\begin{proof}
\ref{PI:AF:stability} is obvious and \ref{PI:AF:sep} is straightforward.
\ref{PI:AF:srepr+sep} and \ref{PI:AF:total-AF} immediately follows from
\ref{PI:AF:sep} and the definitions.
\ref{PI:AF:ample} follows from~\cite[Cor.~4.5.4]{egaII} and \ref{PI:AF:lqf+sep}
follows from~\cite[Thm.~A.2]{laumon} and \ref{PI:AF:ample}.
We deduce \ref{PI:AF:comp} from \ref{PI:AF:total-AF} and as the diagonal of $g$
is locally quasi-finite and separated, we have that \ref{PI:AF:first} follows
from \ref{PI:AF:lqf+sep} and \ref{PI:AF:comp}.
\end{proof}

\begin{warning}
The property \AF{} for a morphism is \emph{not local on the base in the Zariski
  topology}. For a counter-example, let $Y$ be a projective three-fold with two
smooth curves $C$ and $C'$ intersecting transversely at two points $P$ and
$Q$. Then let $X$ be Hironaka's proper non-projective three-fold given by first
blowing-up $C$ and then the strict transform of $C'$ in a neighborhood of $P$
and by first blowing-up $C'$ and then the strict transform of $C$ in a
neighborhood of $Q$. It is well-known that $X$ is not \AF{}, so the morphism
$\map{p}{X}{Y}$ is not \AF{}. On the other hand, $p|_{Y\setminus P}$ and
$p|_{Y\setminus Q}$ are compositions of two blow-ups and hence projective.
\end{warning}

There is an analogue of Chevalley's theorem for \AF{}-schemes: if
$\map{p}{Z}{X}$ is a finite and surjective morphism of algebraic spaces and $Z$
is \AF{}, then so is $X$, cf.\ \cite[Thm.~5.1.5]{gross_thesis} or
Koll\'ar~\cite[Cor.~48]{kollar_finite-equiv-rels}.

The following criterion for projectivity was conjectured by Chevalley and
proved by Kleiman~\cite{kleiman_num_theory_ampleness}.

\begin{theorem}[Chevalley--Kleiman criterion]
Let $X/k$ be a proper regular algebraic scheme. Then $X$ is projective if and
only if $X$ is an \AF{}-scheme.
\end{theorem}

The \AF{} condition is therefore also known as the Chevalley--Kleiman
property~\cite[Def.~47]{kollar_finite-equiv-rels}. Note that proper singular
\AF{}-schemes need not be projective. In fact, there are singular, proper but
non-projective \AF{}-surfaces~\cite{horrocks_nonproj_surface,nori_varieties-without-smooth-emb}.
\end{section}


\begin{section}{Strong homeomorphisms}\label{A:strong-homeomorphisms}

If $\map{f}{X}{Y}$ is a homeomorphism of \emph{topological spaces}, then
the diagonal map is a homeomorphism. If $\map{f}{X}{Y}$ is a universal
homeomorphism of \emph{schemes}, then the diagonal morphism is a universal
homeomorphism. Indeed, it is a surjective immersion, i.e., a nil-immersion.

If $\map{f}{X}{Y}$ is a universal homeomorphism of \emph{algebraic
spaces}, however, then the diagonal is universally bijective but not
necessarily a universal
homeomorphism. A counterexample is given by $Y$ as the affine line and $X$ as
a non-locally separated line. This motivates the following definition.

\begin{definition}
A morphism of algebraic stacks $\map{f}{\stX}{\stY}$ is a
\emph{strong homeomorphism} if $f$ is a universal homeomorphism and
the diagonal $\Delta_f$ is universally submersive.
\end{definition}

\begin{proposition} 
Let $\map{f}{\stX}{\stY}$ and $\map{g}{\stY}{\stZ}$ be morphisms of algebraic
stacks.
\begin{enumerate}
\item A separated universal homeomorphism is a strong homeomorphism.
\item If $f$ is a representable universal homeomorphism, then $f$ is a strong
homeomorphism if and only if $f$ is locally separated, or equivalently, if and
only if $f$ is separated.
\item If $f$ and $g$ are strong homeomorphisms, then so is $g\circ f$.
\item If $g\circ f$ is a strong homeomorphism and
$f$ is universally submersive, then $g$ is a strong homeomorphism.
\item If $g\circ f$ is a strong homeomorphism and
$g$ is a representable strong homeomorphism, then $f$ is a strong
homeomorphism.
\end{enumerate}
\end{proposition}
\begin{proof}
Standard.
\end{proof}

\begin{corollary}\label{C:strong-homeo-and-separatedness}
Let $\stX$ be an algebraic stack and let $X$ and $Y$ be algebraic spaces
together
with morphisms $\map{f}{\stX}{X}$ and $\map{g}{X}{Y}$. If $g\circ f$ is
a strong homeomorphism and $f$ is a universal homeomorphism, then $g$
is separated.
\end{corollary}


\end{section}


\begin{section}{Weil restriction}\label{A:Weil-restrict}
In Section~\ref{S:finite-quotients}, we will have use of the
\emph{Weil restriction},
cf.~\cite[\S7.6]{raynaud-bosch-lutkebohmert_Neron-models}, which is defined
as follows.

\begin{definition}
Let $X\to S$ and $Z\to X$ be morphisms of algebraic spaces. The Weil
restriction of $Z$ along
$X\to S$ is the contravariant functor $\weilr_{X/S}(Z)$ which takes an
$S$-scheme $T$
to the set $\Hom_{X\times_S T}({X\times_S T}, {Z\times_S T})$, that is, the
set of sections of $Z\times_S T\to X\times_S T$.
\end{definition}

\begin{remark}
When the functor $\weilr_{X/S}(Z)$ is representable by a scheme, it is denoted
by
$\Sect_{X/S}Z$ in~\cite[No.~195, \S C 2]{fga}. Let $X\to S$ and $Y\to S$ be
algebraic spaces. Then $\weilr_{X/S}(Z)$ is a generalization of the functor
$T\mapsto \Hom_T(X_T,Y_T)$ where $X_T=X\times_S T$ and $Y_T=Y\times_S T$. In
fact, $\Hom_T(X_T,Y_T)=\weilr_{X/S}(X\times_S Y)(T)$.
\end{remark}

It is easily seen that if $X\to S$ is flat, proper and of finite presentation
and $Z\to X$ is separated,
then $\weilr_{X/S}(Z)$ is an open subfunctor of the Hilbert functor
$\Hilb(Z/S)$. If in addition $Z\to S$ is locally of finite
presentation, then $\weilr_{X/S}(Z)$ is an algebraic space, locally of finite
presentation and separated over $S$. In fact, Artin has shown the algebraicity
of
$\Hilb(Z/S)$ under these hypotheses~\cite[Cor.~6.2]{artin_alg_formal_moduli_I}.
If $X\to S$ is flat, \emph{finite} and of finite presentation and $Z\to X$ is
\emph{arbitrary}, then $\weilr_{X/S}(Z)$ is also an algebraic
space~\cite{rydh_hilbert}.

\begin{proposition}\label{P:PI-prop}
Let $\map{f}{X}{S}$ be a flat and proper morphism of finite presentation between
algebraic spaces and let $Z\to X$ be an \etale{} and separated morphism. Then
$\weilr_{X/S}(Z)\to S$ is \etale{} and separated. If $X\to S$ has a
section, then there is an \etale{} and separated morphism $\weilr_{X/S}(Z)\to
Z\times_X S$.
\end{proposition}
\begin{proof}
As $\weilr_{X/S}(Z)\to S$ is locally of finite presentation and separated, it
is enough to show that $\weilr_{X/S}(Z)\to S$ is also formally \etale{}. This
follows easily from the functorial description of $\weilr_{X/S}(Z)$.
Actually, if we identify \etale{} morphisms with sheaves of sets in the small
\etale{} site, then $\weilr_{X/S}(Z)$ is nothing but the \etale{} sheaf $f_*Z$.
%
If $X\to S$ has a section $\map{s}{S}{X}$ and $T$ is any $S$-scheme, then there
is a natural map
$$\Hom_X(X\times_S T,Z) \to \Hom_X(T,Z)=\Hom_S(T,Z\times_X S)$$
which induces an $S$-morphism $\weilr_{X/S}(Z)\to Z\times_X S$. As
$\weilr_{X/S}(Z)$ and $Z\times_X S$ are \etale{} over $S$, it follows that
$\weilr_{X/S}(Z)\to Z\times_X S$ is \etale{}.
\end{proof}

\end{section}

\bibliography{quotients}

\def\cprime{$'$}
\providecommand{\bysame}{\leavevmode\hbox to3em{\hrulefill}\thinspace}
\providecommand{\MR}{\relax\ifhmode\unskip\space\fi MR }
\providecommand{\MRhref}[2]{%
  \href{http://www.ams.org/mathscinet-getitem?mr=#1}{#2}
}
\providecommand{\href}[2]{#2}
\begin{thebibliography}{EGA\textsubscript{IV}}

\bibitem[AOV08]{abramovich_olsson_vistoli_tame_stacks}
Dan Abramovich, Martin Olsson, and Angelo Vistoli, \emph{Tame stacks in
  positive characteristic}, Ann. Inst. Fourier (Grenoble) \textbf{58} (2008),
  no.~4, 1057--1091,
  \href{http://arXiv.org/abs/math/0703310}{\mbox{arXiv:math/0703310}}.

\bibitem[Art69a]{artin_alg_formal_moduli_I}
M.~Artin, \emph{Algebraization of formal moduli. {I}}, Global Analysis (Papers
  in Honor of K. Kodaira), Univ. Tokyo Press, Tokyo, 1969, pp.~21--71.

\bibitem[Art69b]{artin_implicit_function_thm}
\bysame, \emph{The implicit function theorem in algebraic geometry}, Algebraic
  Geometry (Internat. Colloq., Tata Inst. Fund. Res., Bombay, 1968), Oxford
  Univ. Press, London, 1969, pp.~13--34.

\bibitem[Art71]{artin_Hensel-join}
\bysame, \emph{On the joins of {H}ensel rings}, Advances in Math. \textbf{7}
  (1971), 282--296 (1971).

\bibitem[Art73]{artin_theoremes_de_repr}
\bysame, \emph{Th\'eor\`emes de repr\'esentabilit\'e pour les espaces
  alg\'ebriques}, Les Presses de l'Universit\'e de Montr\'eal, Montreal, Que.,
  1973, En collaboration avec Alexandru Lascu et Jean-Fran\c cois Boutot,
  S\'eminaire de Math\'ematiques Sup\'erieures, No. 44 (\'Et\'e, 1970).

\bibitem[Art74]{artin_versal_def_alg_stacks}
\bysame, \emph{Versal deformations and algebraic stacks}, Invent. Math.
  \textbf{27} (1974), 165--189.

\bibitem[BLR90]{raynaud-bosch-lutkebohmert_Neron-models}
Siegfried Bosch, Werner L{\"u}tkebohmert, and Michel Raynaud, \emph{N\'eron
  models}, Springer-Verlag, Berlin, 1990.

\bibitem[Con05]{conrad_coarsespace}
Brian Conrad, \emph{The {K}eel--{M}ori theorem via stacks}, Nov 2005, Preprint,
  p.~12.

\bibitem[DG70]{gabriel_algebraic_groups}
Michel Demazure and Pierre Gabriel, \emph{Groupes alg\'ebriques. {T}ome {I}:
  {G}\'eom\'etrie alg\'ebrique, g\'en\'eralit\'es, groupes commutatifs}, Masson
  \& Cie, \'Editeur, Paris, 1970, Avec un appendice {\it Corps de classes
  local}\ par Michiel Hazewinkel.

\bibitem[EGA\textsubscript{I}]{egaI_NE}
A.~Grothendieck, \emph{\'{E}l\'ements de g\'eom\'etrie alg\'ebrique. {I}. {L}e
  langage des sch\'emas}, second ed., Die Grundlehren der mathematischen
  Wissenschaften in Einzeldarstellungen, vol. 166, Springer-Verlag, Berlin,
  1971.

\bibitem[EGA\textsubscript{II}]{egaII}
\bysame, \emph{\'{E}l\'ements de g\'eom\'etrie alg\'ebrique. {II}. \'{E}tude
  globale \'el\'ementaire de quelques classes de morphismes}, Inst. Hautes
  \'Etudes Sci. Publ. Math. (1961), no.~8, 222.

\bibitem[EGA\textsubscript{IV}]{egaIV}
\bysame, \emph{\'{E}l\'ements de g\'eom\'etrie alg\'ebrique. {IV}. \'{E}tude
  locale des sch\'emas et des morphismes de sch\'emas}, Inst. Hautes \'Etudes
  Sci. Publ. Math. (1964-67), nos.~20, 24, 28, 32.

\bibitem[ES04]{skjelnes_ekedahl_good_component}
Torsten Ekedahl and Roy Skjelnes, \emph{Recovering the good component of the
  {H}ilbert scheme}, Preprint, May 2004,
  \href{http://arXiv.org/abs/math.AG/0405073}{\mbox{arXiv:math.AG/0405073}}.

\bibitem[Fer03]{ferrand_conducteur}
Daniel Ferrand, \emph{Conducteur, descente et pincement}, Bull. Soc. Math.
  France \textbf{131} (2003), no.~4, 553--585.

\bibitem[FGA]{fga}
A.~Grothendieck, \emph{Fondements de la g\'eom\'etrie alg\'ebrique. [{E}xtraits
  du {S}\'eminaire {B}ourbaki, 1957--1962.]}, Secr\'etariat math\'ematique,
  Paris, 1962.

\bibitem[Gab63]{gabriel_quotients}
Pierre Gabriel, \emph{Construction de pr\'esch\'emas quotient}, Sch\'emas en
  Groupes (S\'em. G\'eom\'etrie Alg\'ebrique, Inst. Hautes \'Etudes Sci.,
  1963/64), Fasc. 2a, Expos\'e 5, Inst. Hautes \'Etudes Sci., Paris, 1963,
  p.~37.

\bibitem[Gir64]{giraud_descente}
Jean Giraud, \emph{M\'ethode de la descente}, Bull. Soc. Math. France M\'em.
  \textbf{2} (1964), viii+150.

\bibitem[GIT]{GIT}
D.~Mumford, J.~Fogarty, and F.~Kirwan, \emph{Geometric invariant theory}, third
  ed., Springer-Verlag, Berlin, 1994.

\bibitem[GM11]{geer-moonen_abel-vars}
Gerard van~der Geer and Ben Moonen, \emph{{A}belian {V}arieties}, Book in
  preparation, 2011.

\bibitem[Gro10]{gross_thesis}
Philipp Gross, \emph{Vector bundles as generators on schemes and stacks},
  Ph{D}. {T}hesis, D\"usseldorf, May 2010.

\bibitem[Hor71]{horrocks_nonproj_surface}
G.~Horrocks, \emph{Birationally ruled surfaces without embeddings in regular
  schemes}, Bull. London Math. Soc. \textbf{3} (1971), 57--60.

\bibitem[Kle66]{kleiman_num_theory_ampleness}
Steven~L. Kleiman, \emph{Toward a numerical theory of ampleness}, Ann. of Math.
  (2) \textbf{84} (1966), 293--344.

\bibitem[KM97]{keel_mori_quotients}
Se{\'a}n Keel and Shigefumi Mori, \emph{Quotients by groupoids}, Ann. of Math.
  (2) \textbf{145} (1997), no.~1, 193--213,
  \href{http://arXiv.org/abs/alg-geom/9508012}{\mbox{arXiv:alg-geom/9508012}}.

\bibitem[Knu71]{knutson_alg_spaces}
Donald Knutson, \emph{Algebraic spaces}, Springer-Verlag, Berlin, 1971, Lecture
  Notes in Mathematics, Vol. 203.

\bibitem[Kol97]{kollar_quotients}
J{\'a}nos Koll{\'a}r, \emph{Quotient spaces modulo algebraic groups}, Ann. of
  Math. (2) \textbf{145} (1997), no.~1, 33--79,
  \href{http://arXiv.org/abs/alg-geom/9503007}{\mbox{arXiv:alg-geom/9503007}}.

\bibitem[Kol11]{kollar_finite-equiv-rels}
J\'anos Koll\'ar, \emph{Quotients by finite equivalence relations}, Current
  {D}evelopments in {A}lgebraic {G}eometry, Math. Sci. Res. Inst. Publ.,
  vol.~59, Cambridge Univ. Press, Cambridge, 2011, pp.~227--256.

\bibitem[LMB00]{laumon}
G{\'e}rard Laumon and Laurent Moret-Bailly, \emph{Champs alg\'ebriques},
  Springer-Verlag, Berlin, 2000.

\bibitem[Mat76]{matsuura_quotients}
Yutaka Matsuura, \emph{On a construction of quotient spaces of algebraic
  spaces}, Proceedings of the Institute of Natural Sciences, Nihon University
  \textbf{11} (1976), 1--6.

\bibitem[Mil80]{milne_etale_coh}
James~S. Milne, \emph{\'{E}tale cohomology}, Princeton Mathematical Series,
  vol.~33, Princeton University Press, Princeton, N.J., 1980.

\bibitem[Mum70]{mumford_abel_vars}
David Mumford, \emph{Abelian varieties}, Tata Institute of Fundamental Research
  Studies in Mathematics, No. 5, Published for the Tata Institute of
  Fundamental Research, Bombay, 1970.

\bibitem[Nor78]{nori_varieties-without-smooth-emb}
M.~V. Nori, \emph{Varieties with no smooth embeddings}, C. {P}. {R}amanujam---a
  tribute, Tata Inst. Fund. Res. Studies in Math., vol.~8, Springer, Berlin,
  1978, pp.~241--246.

\bibitem[Ray70]{raynaud_hensel_rings}
Michel Raynaud, \emph{Anneaux locaux hens\'eliens}, Lecture Notes in
  Mathematics, Vol. 169, Springer-Verlag, Berlin, 1970.

\bibitem[RS10]{rydh-skjelnes_gen_etale_fam}
David Rydh and Roy Skjelnes, \emph{An intrinsic construction of the principal
  component of the {H}ilbert scheme}, J. Lond. Math. Soc. (2) \textbf{82}
  (2010), no.~2, 459--481.

\bibitem[Ryd08a]{rydh_thesis}
David Rydh, \emph{Families of cycles and the {C}how scheme}, Ph.D. thesis,
  Royal Institute of Technology, Stockholm, May 2008, p.~218.

\bibitem[Ryd08b]{rydh_famzerocycles-I}
\bysame, \emph{Families of zero-cycles and divided powers: {I}.
  {R}epresentability}, Preprint, Part of~\cite{rydh_thesis}, Mar 2008,
  \href{http://arXiv.org/abs/0803.0618v1}{\mbox{arXiv:0803.0618v1}}.

\bibitem[Ryd08c]{rydh_gammasymchow_inprep}
\bysame, \emph{{H}ilbert and {C}how schemes of points, symmetric products and
  divided powers}, Part of~\cite{rydh_thesis}, May 2008.

\bibitem[Ryd09]{rydh_noetherian-approx}
\bysame, \emph{Noetherian approximation of algebraic spaces and stacks},
  Preprint, Apr 2009,
  \href{http://arXiv.org/abs/0904.0227v2}{\mbox{arXiv:0904.0227v2}}.

\bibitem[Ryd10]{rydh_submersion_and_descent}
\bysame, \emph{Submersions and effective descent of \'etale morphisms}, Bull.
  Soc. Math. France \textbf{138} (2010), no.~2, 181--230.

\bibitem[Ryd11a]{rydh_etale-devissage}
\bysame, \emph{\'{E}tale d\'evissage, descent and pushouts of stacks}, J.
  Algebra \textbf{331} (2011), 194--223,
  \href{http://arXiv.org/abs/1005.2171}{\mbox{arXiv:1005.2171}}.

\bibitem[Ryd11b]{rydh_hilbert}
\bysame, \emph{Representability of {H}ilbert schemes and {H}ilbert stacks of
  points}, Comm. Algebra \textbf{39} (2011), no.~7, 2632--2646.

\bibitem[Sch03]{schroer_bigger-Brauer}
Stefan Schr{\"o}er, \emph{The bigger {B}rauer group is really big}, J. Algebra
  \textbf{262} (2003), no.~1, 210--225.

\bibitem[Ser59]{serre_groupes_alg_corps_de_classes}
Jean-Pierre Serre, \emph{Groupes alg\'ebriques et corps de classes},
  Publications de l'institut de math\'ematique de l'universit\'e de Nancago,
  VII. Hermann, Paris, 1959.

\bibitem[SGA\textsubscript{1}]{sga1}
A.~Grothendieck (ed.), \emph{Rev\^etements \'etales et groupe fondamental},
  Springer-Verlag, Berlin, 1971, S\'eminaire de G\'eom\'etrie Alg\'ebrique du
  Bois Marie 1960--1961 (SGA 1), Dirig\'e par Alexandre Grothendieck.
  Augment\'e de deux expos\'es de M. Raynaud, Lecture Notes in Mathematics,
  Vol. 224.

\bibitem[SGA\textsubscript{4}]{sga4}
M.~Artin, A.~Grothendieck, and J.~L. Verdier (eds.), \emph{Th\'eorie des topos
  et cohomologie \'etale des sch\'emas}, Springer-Verlag, Berlin, 1972--1973,
  S\'eminaire de G\'eom\'etrie Alg\'ebrique du Bois Marie 1963--1964 (SGA 4).
  Dirig\'e par M. Artin, A. Grothendieck et J. L. Verdier. Avec la
  collaboration de P. Deligne et B. Saint-Donat, Lecture Notes in Mathematics,
  Vol. 269, 270, 305.

\bibitem[SP]{stacks-project}
{The Stacks Project Authors}, \emph{Stacks project},
  \url{http://math.columbia.edu/algebraic_geometry/stacks-git}.

\bibitem[Vis05]{vistoli_descent}
Angelo Vistoli, \emph{Grothendieck topologies, fibered categories and descent
  theory}, Fundamental algebraic geometry, Math. Surveys Monogr., vol. 123,
  Amer. Math. Soc., Providence, RI, 2005, pp.~1--104.

\bibitem[W{\l}o93]{wlodarczyk_toric-embeddings}
Jaros{\l}aw W{\l}odarczyk, \emph{Embeddings in toric varieties and
  prevarieties}, J. Algebraic Geom. \textbf{2} (1993), no.~4, 705--726.

\end{thebibliography}
\bibliographystyle{dary}

\end{document}